# CLASSICAL TWO-PARABOLIC $T$-SCHOTTKY GROUPS

*By*

JANE GILMAN*AND PETER WATERMAN

*Dedicated to the memory of Carol Feltz Waterman*

**Abstract.** A $T$-Schottky group is a discrete group of Möbius transformations whose generators identify pairs of possibly-tangent Jordan curves on the complex sphere $\hat{\mathbb{C}}$. If the curves are Euclidean circles, then the group is termed classical $T$-Schottky.

We describe the boundary of the space of classical $T$-Schottky groups affording two parabolic generators within the larger parameter space of all $T$-Schottky groups with two parabolic generators. This boundary is surprisingly different from that of the larger space. It is analytic, while the boundary of the larger space appears to be fractal. Approaching the boundary of the smaller space does not correspond to pinching; circles necessarily become tangent, but extra parabolics need not develop.

As an application, we construct an explicit one parameter family of two parabolic generator non-classical $T$-Schottky groups.

## 1 Introduction

We are interested in two types of Möbius groups that are defined by the geometry of their actions on $\hat{\mathbb{C}}$: $T$-Schottky groups and nsdc groups.

After a general overview, in Section 4 and later, we restrict our attention to groups possessing two parabolic generators, unless explicitly stated otherwise.

### 1.1 $T$-Schottky groups

**Definition 1.1.** Let $C_i, C'_i, i = 1, \ldots, n$ be a set of $2n$ circles in $\hat{\mathbb{C}}$ such that the interiors of the $2n$ circles are all pairwise disjoint. Let $\mathcal{F}$ be the intersection of the exteriors of these circles.

- $\mathcal{F}$ is called a **classical $T$-Schottky domain**.

*Research supported in part by NSA grant #02G-186, the NSF, IHES and the Rutgers Research Council.

- A **classical $T$-Schottky pairing** is a set of $n$ Möbius transformations, $g_1, \ldots, g_n$, the **side pairings**, where each $g_i$ maps the exterior of $C_i$ onto the interior of $C_i'$.
- The set of circles together with the side pairings is called a **classical $T$-Schottky configuration**.
- The group generated by the side-pairing is called a **marked classical $T$-Schottky group** on the generators $g_1, \ldots, g_n$.
- A group of Möbius transformations is a **classical $T$-Schottky group** if it is a classical marked $T$-Schottky group on *some* set of generators.

More generally, one can relax the requirement that the sides be Euclidean circles.

**Definition 1.2.**

- A **marked $T$-Schottky group** is a marked classical $T$-Schottky group except that the $2n$ curves are allowed to be Jordan curves with at most finitely many tangencies.
- A group of Möbius transformations is a $T$**-Schottky group** if it is a marked $T$-Schottky group on *some* set of generators.
- A marked $T$-Schottky group is **marked non-classical** if it is not classical on the marked set of generators.
- A $T$-Schottky group is **non-classical** if it is not classical on any set of generators.

Recall that a marking for a group is a choice of an ordered set of generators for the group. The groups we are primarily interested in are two generator groups, and it will be clear from the context whether we are thinking of the group generated by the given set of generators or the marked group with the given set of generators as the marking.

We note that corresponding to any given marking, there typically will be many different $T$-Schottky domains, if there is one. Further, a given classical $T$-Schottky group will typically have many possible $T$-Schottky configurations corresponding to different markings. We further point out that one may be able to associate different side pairings to a given set of classical $T$-Schottky circles. Moreover, one can give a non-classical pairing for a classical $T$-Schottky group. Indeed, it is also possible for a classical $T$-Schottky group to be non-classical on some marked set of generators, as in Zarrow [33]. The effect of changing generators is discussed in Section 6.1.

We remind the reader that if the sides of a $T$-Schottky domain are disjoint, then it is well-known that the resulting group is purely loxodromic. This follows easily by noting that, up to conjugacy, any transformation has one fixed point inside or on one of the Jordan curves and its inverse a fixed point inside or on *another*. Groups possessing a domain for which the curves are disjoint are termed **Schottky** .

Groups with parabolic elements possessing a $T$-Schottky domain with tangencies *solely at parabolic fixed points lying on paired circles* are termed of **Schottky Type** or **Noded Schottky**.

Maskit ([24] page 314) has shown that a finitely generated, free, purely loxodromic Kleinian group (of the second kind) is a Schottky group. Similarly, a geometrically finite, free Kleinian group (of the second kind) with parabolics is of Schottky Type. As a consequence, a $T$-Schottky group is either a Schottky group or a Schottky Type group. However, a classical $T$-Schottky group could fail to be either a classical Schottky group or a classical Schottky Type group: it might not be possible to remove the non-parabolic tangencies and keep circles for the sides. This is one of the phenomena we examine in this paper.

The existence of non-classical Schottky groups was proved by Marden [20] and an explicit example constructed by Yamamoto [32].

Jørgensen, Marden and Maskit [14] have shown that groups on the boundary of classical Schottky space are geometrically finite. Further, in [26], Maskit has shown that classical Schottky space is connected.

Maskit and Swarup [27] show that any torsion-free two parabolic generator Kleinian group of the second kind is geometrically finite. It is these groups that we study in this paper.

Once general tangencies are permitted, it is possible that a $T$-Schottky domain $\mathcal{F}$ will have a loxodromic fixed point on its boundary and will then fail to be a fundamental domain. Nevertheless, no points of $\mathcal{F}$ are group equivalent; so existence of a $T$-Schottky domain guarantees discreteness with non-empty ordinary set. However, it turns out that in spite of having non-parabolic tangencies, the *extreme* domains we consider are, in fact, fundamental domains.

Our main results give an explicit description as to when such a group is classical $T$-Schottky (Theorem 4.2, 6.18) and thereby enable us to construct an explicit one parameter family of two parabolic generator non-classical $T$-Schottky groups (Theorem 7.2). Conjugacy classes of two parabolic generator groups are parametrized by the one complex parameter $\lambda$, where $4\lambda^2 + 2$ is the trace of the commutator of the Möbius transformations that generate the group. We describe classical $T$-Schottky space by describing this group invariant, that is, by describing those $\lambda$ that lie on the boundary of the space. Specifically, if such a $\lambda = x + iy$, then $|y| = 1 - x^2/4$.

**1.2  Non-separating disjoint circle groups**   There is then a natural variation on the concept of classical $T$-Schottky group, due to Gilman [8], which can be easier to utilize: **non-separating disjoint circle groups**.

**1.2.1  Notation and terminology**   An element of $PSL(2,\mathbb{C})$ acts as a Möbius transformation on $\widehat{\mathbb{C}}$, the complex sphere. The action extends in a natural way to upper-half three-space $\mathbb{H} = \{(x, y, t) : t > 0, \quad x, y, t, \in \mathbb{R}\}$ and the extended transformation is an isometry acting on $\mathbb{H}$ endowed with the hyperbolic metric. When considered as the boundary of hyperbolic three-space, $\widehat{\mathbb{C}}$ is referred to as the sphere at infinity. Elements of $PSL(2,\mathbb{C})$ are classified as loxodromic, elliptic, parabolic according to the square of their traces. The classification of transformations can also be described by the action on $\widehat{\mathbb{C}}$ or on $\widehat{\mathbb{C}} \cup \mathbb{H}$. We follow the notation of [6]. For $x$ and $y \in (\widehat{\mathbb{C}} \cup \mathbb{H})$ with $x \neq y$, we let $[x, y]$ denote the oriented hyperbolic geodesic passing through $x$ and $y$. The *ends* of the geodesic $[x, y]$ are by definition the points $v$ and $w$ with $\{v, w\} = [x, y] \cap \widehat{\mathbb{C}}$. If $v$ and $w$ are the ends of $[x, y]$, then $[x, y] = [v, w]$. The notation $[x, y]$ also indicates a direction, so that $[y, x]$ is the geodesic with the opposite orientation.

Following Fenchel [6], we include *improper lines* in our considerations. A *proper line* in one whose ends are distinct. An improper line is one whose ends coincide.

A Möbius transformation fixes one or two points on $\widehat{\mathbb{C}}$. A loxodromic transformation $A$ has two distinct fixed points, as does an elliptic transformation. These transformations fix the geodesic $X$ in $\mathbb{H}$ whose ends are the fixed points of $A$ on $\widehat{\mathbb{C}}$. This geodesic is called the *axis* of $A$ and is denoted by $Ax_A$. A parabolic transformation has one fixed point. Using the terminology of improper lines, parabolic transformations also have axes. If $A$ is a parabolic transformation with fixed point $u \in \widehat{\mathbb{C}}$, we consider the improper line $[u, u]$ to be its axis.

A subgroup $G$ of $PSL(2,\mathbb{C})$ is termed *elementary* if there is a finite G-orbit in $\widehat{\mathbb{H}}$ ([29]). In this paper, we consider only non-elementary groups. According to this definition, a group generated by a loxodromic and another element which share exactly one fixed point is elementary. It is well known that such a group is not discrete.

With this terminology, every pair of distinct lines (proper or improper) in $\mathbb{H}$ has a unique common perpendicular (proper or improper) (III.2.6 of [6]). Further, every pair of elements $A, B \in PSL(2,\mathbb{C})$ with distinct axes determines a unique (proper or improper) hyperbolic line in $\mathbb{H}$, the common perpendicular of their axes; we denote this common perpendicular by $L$ and the ends of $L$ by $n$ and $n'$. Note that $L$ is improper if $n = n'$.

**1.2.2  Half-turns**  For any hyperbolic line $[x, y]$ with $x \neq y$, we let $H_{[x,y]}$ be the half-turn about the line with ends $x$ and $y$. We note that $H_{[x,y]}$ leaves invariant every hyperbolic plane $\mathbb{P}$ whose boundary $C_\mathbb{P}$, also called its *horizon*, passes through $x$ and $y$ ([8]). The half-turn interchanges the exterior and the interior of $\mathbb{P} \in \mathbb{H}$.

**1.2.3  Definition**  When the common perpendicular to the axes of $A$ and $B$ is a proper line, there is a natural construction that associates an ordered six-tuple of complex numbers $(a, a', n, n', b, b')$ to the ordered pair of transformations. Conversely, an ordered six-tuple in $\widehat{\mathbb{C}}$ $(a, a', n, n', b, b')$ with $a \neq a', b \neq b', n \neq n'$ determines an ordered pair of transformations.

Let $G = \langle A, B \rangle$ be a non-elementary group, so that $L = [n, n']$, the common perpendicular to the axes of $A$ and $B$, is a proper line. There are unique hyperbolic lines $L_A$ and $L_B$ such that $A = H_{L_A} \cdot H_L$ and $B = H_{L_B} \cdot H_L$. We let $a$ and $a'$ be the ends of $L_A$, so that $L_A = [a, a']$, and let $b$ and $b'$ be those of $L_B$, so that $L_B = [b, b']$.

**Definition 1.3.** The marked *three generator group* determined by $G = \langle A, B \rangle$ is denoted by $\mathbb{T}G$ and defined by $\mathbb{T}G = \langle H_{L_A}, H_L, H_{L_B} \rangle$.

By construction, $G$ is a normal subgroup of $\mathbb{T}G$ of index at most two, which immediately implies that $G$ is discrete if and only if $\mathbb{T}G$ is discrete.

We also define

**Definition 1.4.** The **ortho-end** of the marked non-elementary group $G = \langle A, B \rangle$ is the six-tuple of complex numbers $(a, a', n, n', b, b')$.

**Definition 1.5.** Six points in $(a, a', n, n', b, b') \in \widehat{\mathbb{C}}^6$ such that $a \neq a'$, $b \neq b'$ and $n \neq n'$ have the **non-separating disjoint circle property** if there are pairwise disjoint or tangent circles on $\widehat{\mathbb{C}}$, $C_A$, $C_D$ and $C_B$ (respectively) passing through $a$ and $a'$, $n$ and $n'$, and $b$ and $b'$ (respectively) such that no one circle separates the other two.

**Definition 1.6.** A non-elementary marked group $G = \langle A, B \rangle$ is a marked *non-separating disjoint circle group*, a marked **nsdc** group, for short, if the ortho-end of $A$ and $B$ has the non-separating disjoint circle property. Further $G$ is a *non-separating disjoint circle group* (an *nsdc* group for short) if some pair of generators for $G$ has the non-separating disjoint circle property.

**Remark 1.7.** We note that if one of $A$ or $B$ is elliptic, then the ortho-end of the marked group is well-defined but does not have the nsdc property. Note that the condition $n \neq n'$ in the definition of an nsdc group rules out the possibility of $A$ and $B$ sharing a fixed point. Such groups do not determine any three generator group and do not have an ortho-end and thus are not nsdc. Note that if $n = n'$, in the case one of the generators is loxodromic and the two generators share a fixed point, the group is not discrete; if one is elliptic, the group is not nsdc; and if both are parabolic with distinct axes, $n \neq n'$.

If $G$ is nsdc, we also call the corresponding group $\mathbb{T}G$ an nsdc group.

Motivating the formulation of the nsdc concept was the fact that the nsdc property implies that a group is discrete. Thus it may be seen as a discreteness test. Since each circle is mapped to itself under the corresponding half-turn, existence of circles with the nsdc property guarantees that $\mathbb{T}G$ is discrete and a free product of the three half-turns, (although, as in the $T$-Schottky case, the exterior of the circles may fail to be a fundamental domain). Hence, a marked nsdc group is always discrete and free.

Note that if $G = \langle A, B \rangle$ is marked nsdc, then it is marked classical $T$-Schottky [8] since $A$ pairs the circles $C_A$ and $H_N(C_A)$ and $B$ the circles $C_B$ and $H_N(C_B)$. If $G$ is nsdc, then $\mathbb{T}G$ is a degree two extension.

## 2  Method

Our method is to analyze a geometric configuration of circles on the boundary of $\mathbb{H}$. This is appropriate because it is by the existence of such circles that $T$-Schottky groups and non-separating disjoint circle groups are defined.

We find geometric configurations that have extra tangencies, that is, those one might expect to be configurations that give a boundary group. However, we prove that not all of these potential boundary groups are actual boundary groups: they have neighborhoods where the tangencies can be pulled apart.

For example, given the fundamental domain (Figure 1) for $\langle z \mapsto -2z \rangle$, we may remove the tangency by replacing the circle $C_1$ and its image $C_2$ with the Jordan curve $J_1$ and its image $J_2$ shown in Figure 2. The Jordan curve $J_1$ is chosen to agree with the circle $C_1$ except near the point of tangency. Alternatively, one can replace $C_1$ with a nearby circle $C_1'$ and its image circle $C_2'$. We are interested in this second situation. It is shown in Figure 3. As the number of sides and tangencies increases, the analysis becomes more subtle.

Our main tool is simply to relate the geometric configurations to algebraic parameters and to use some analytic tools to find boundaries.

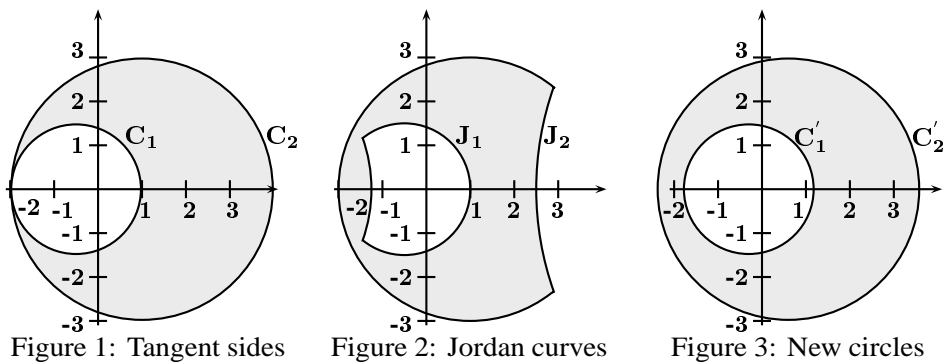

Figure 1: Tangent sides     Figure 2: Jordan curves     Figure 3: New circles

## 3 Preliminaries

**3.1 Möbius Transformations**   A $2 \times 2$ non-singular matrix $M = \begin{pmatrix} a & b \\ c & d \end{pmatrix}$ with complex entries $a$, $b$, $c$, $d$ acts on $\widehat{\mathbb{C}}$, *the sphere at infinity*, as a fractional linear transformation. We often do not distinguish notationally between a matrix and the mapping it induces on $\widehat{\mathbb{C}}$.

Since a matrix and its negative determine the same action, all fractional linear transformations can be considered as elements of $PSL(2, \mathbb{C})$. Unless otherwise required by the situation, we assume pull backs of elements to $SL(2, \mathbb{C})$ have positive trace. In particular, for a two generator group, once the signs of the trace of the (pull backs) of the generators are chosen, the signs of the races of all other elements of the group are determined.

Elements of $SL(2, \mathbb{C})$, or equivalently $PSL(2, \mathbb{C})$, are classified by the square of their traces. We use the term *loxodromic* to include both purely hyperbolic and strictly loxodromic transformations. We let $\text{tr}^2(A)$ denote $(\text{tr}(A))^2$; the trace of $M$ is $(a + d)/\sqrt{D}$, where $D = ad - bc \neq 0$.

A parabolic element has a unique fixed point on the boundary of $\mathbb{H}$ which will often be referred to as a **cusp**.

**3.2 Kleinian groups**   If a subgroup $G$ of $PSL(2, \mathbb{C})$ is discrete, it is called a *Kleinian group*. We term $z \in \widehat{\mathbb{C}}$ an *ordinary point* if it has a neighborhood $U$ such that $g(U) \cap U \neq \emptyset$ for at most finitely many $g \in G$. The set of ordinary points is denoted by $\Omega(G)$ and is also known as the *regular set* or *the set of discontinuity*.

A Kleinian group $G$ is said to be of the second kind if $\Omega(G)$ is non-empty. By definition, the limit set of $G$ is the complement in $\widehat{\mathbb{C}}$ of the ordinary set and is denoted by $\Lambda(G)$. If $G$ is a non-elementary Kleinian group, then $\Lambda(G)$ is the set of accumulation points of the orbit of any point in $\widehat{\mathbb{C}}$ under $G$.

If $G$ is discrete, we can form the quotient $\mathbb{H}/G$. Provided $G$ contains no elements of finite order (that is, no elliptic elements), it is a three-manifold; otherwise it is an *orbifold*. The conformal boundary of $\mathbb{H}/G$ is determined by the action of $G$ on the sphere at infinity and is the quotient of $\Omega(G)$ under the action of $G$.

It is a consequence of a theorem of Ahlfors [1] that if $G$ is a finitely generated Kleinian group of the second kind, then $\Omega(G)/G$ is an *analytically finite* Riemann surface. That means that $\Omega(G)/G$ is a Riemann surface of finite type or a disjoint union of a finite number of Riemann surfaces of finite type where every boundary component has a neighborhood conformally equivalent to a punctured disc. Even if elliptics are present, the quotient will be a Riemann surface (see Theorem 6.2.1 of [2]). Further, it is shown in [21] that a finitely generated free Kleinian group is geometrically finite if and only if every parabolic fixed point is *doubly cusped*, that is, the punctured discs occur in pairs.

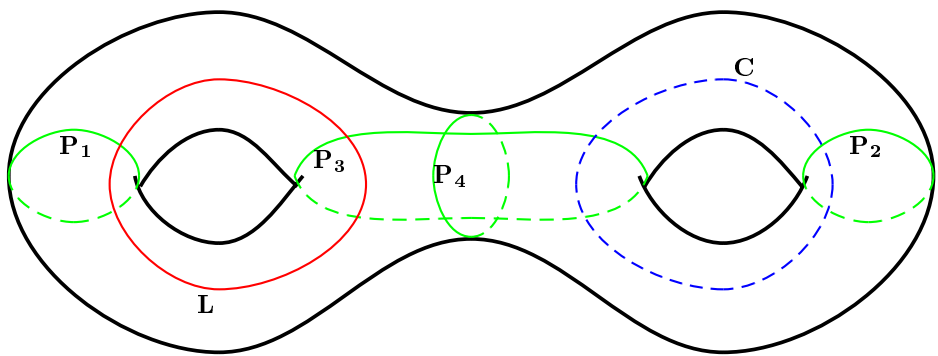

Figure 4: **Some pinchable curves on a surface of genus two.** The surface is the quotient of a rank two Schottky group. The curves $L$ and $C$ are the images of the sides of the Schottky domain. $P_1$ and $P_2$ are curves that can be pinched to form a noded surface of genus two with two double cusps. If $P_4$ is pinched, a pair of once-punctured tori is obtained.

**3.3 Noded surfaces and pinching** A *node* on a Riemann surface is a point with a neighborhood isomorphic to

$$\{z \in \mathbb{C}^2 : |z_1| < 1, |z_2| < 1, z_1 z_2 = 0\}$$

([3, 4]). In particular, a doubly cusped parabolic fixed point corresponds to a node on $\Omega(G)/G$.

If $\Sigma$ is a Riemann surface with nodes, then $\Sigma - \{\text{nodes}\}$ is a disjoint union of Riemann surfaces, called the *parts* of $\Sigma$. We require that every part $S_j$ be a Riemann surface of type $(p_j, n_j)$, where $p_j$ is the genus and $n_j$ the number of punctures and $3p_j - 3 + n_j \geq 0$. The genus of a noded surface is the genus $g$ obtained by thickening the nodes, that is, $g = \sum_j (p_j - 1) + (1/2) \sum_j n_j + 1$. If the noded surface has genus two, the possibilities for number of parts, the genus of the parts and the number of nodes are easily enumerated.

Intuitively, a noded surface can be viewed as obtained from the original surface by *pinching* along certain curves.

To make this notion more precise, recall that if $F$ is any free group, then an element $W$ in $F$ is *primitive* if there is a minimal generating set for $F$ containing it (see [18]). An element of a group is *primary* if it is not a power of any other element of the group.

**Definition 3.1.** Pinching: A surface $\widehat{\Sigma} = \Omega(\widehat{F})/\widehat{F}$ (respectively the group $\widehat{F}$) is obtained from $\Sigma = \Omega(F)/F$ (respectively from $F$) by **pinching** $W$ if there is an isomorphism of $F$ onto $\widehat{F}$ which sends a non-parabolic primary element $W$ to a parabolic.

A result of Yamamoto ([31]), stated for two generator groups, says the following.

Let $F$ be a Schottky group of rank two with loxodromic generators $A_1, A_2$. Consider simple disjoint geodesics $L_1, ..., L_q$, $q \leq 3$, on $\Sigma = \Omega(F)/F$ defined by the words $W_1, ..., W_q$, where each $W_i$ is primary in $F$ and for $i \neq j$ the cyclic subgroups generated by $W_i$ and $W_j$ are not conjugate in $F$. There is then a Schottky type group $\widehat{F}$ and an isomorphism $\psi : F \to \widehat{F}$ where $\psi(W_1), ..., \psi(W_q)$ and their powers and conjugates are exactly the parabolic elements of $\widehat{F}$. Such a set of geodesics $L_1, ..., L_q$ is known as a set of *pinchable geodesics* [10]; we also call the words $W_1, ... W_q$ *pinchable*.

We extend this to groups that have parabolic generators. A group $G^{\mathcal{P}}$ generated by two parabolics comes from a group $G$, with loxodromic generators, by pinching along some curves, say $W_1$ and $W_2$. If $\widehat{G}$ is any group that comes from $G^{\mathcal{P}}$ by further pinching along $W_3$, then $\widehat{G}$ can be obtained from $G$ by pinching along $W_1, W_2$, and $W_3'$ where the curve $W_3'$ is sent to $W_3$ under the isomorphism of $G$ onto $G^{\mathcal{P}}$.

A pinchable word may or may not be primitive. For example, a noded surface with two parts, one a punctured torus and the other a thrice punctured sphere, does not correspond to a group with two parabolic generators. That is, it does not arise from pinching along two primitive curves.

Since $G^{\mathcal{P}}$ is generated by two parabolics, the surface $\Omega(G^{\mathcal{P}})/G^{\mathcal{P}}$ represents a sphere with four punctures. Additional pinching along $W_3$ will yield two thrice punctured spheres when $W_3$ is primitive and two once punctured tori when $W_3$ is not primitive.

## 4  Two parabolic generator groups

From now on, we restrict our attention to non-elementary groups that can be generated by two parabolics. Such groups are determined up to conjugacy by a single complex parameter.

**4.0.1  The Lyndon–Ullman parameter: $\lambda$.**  If $G$ is a marked group with two parabolic generators, $S$ and $T$, then up to $PSL(2, \mathbb{C})$ conjugation, $G = G_\lambda = \langle S, T \rangle$, where $S = \begin{pmatrix} 1 & 0 \\ 1 & 1 \end{pmatrix}$ and $T = T_\lambda = \begin{pmatrix} 1 & 2\lambda \\ 0 & 1 \end{pmatrix}$ for some $\lambda = |\lambda|e^{i\omega} \in \mathbb{C}$. We assume that $\lambda \neq 0$, so that $G$ is non-elementary.

This normalization was used by Lyndon and Ullman [17], who found various regions with $G_\lambda$ free as shown in Figure 5.

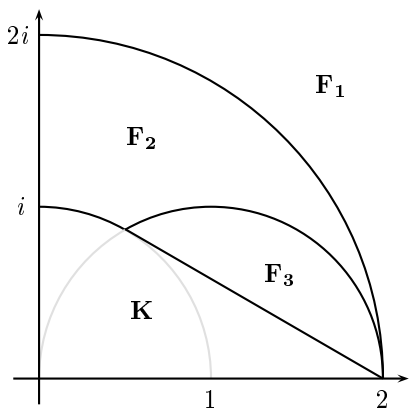

Figure 5: The Lyndon-Ullman Free Region(s) in $\mathbb{C}$.

We note for future reference that since $\operatorname{tr} S = 2, \operatorname{tr} T = 2$ and $\lambda \neq 0$,

(1) $$\operatorname{tr}[S, T] - 2 = \operatorname{tr} STS^{-1}T^{-1} - 2 = 4\lambda^2;$$

(2) $$\operatorname{tr} ST^{-1} = \operatorname{tr} TS^{-1} = 2 - 2\lambda \quad \text{and} \quad \operatorname{tr} ST = \operatorname{tr}(TS)^{-1} = 2 + 2\lambda.$$

(3) $\operatorname{tr} ST = \pm 2 \Leftrightarrow \lambda = -2, \ \operatorname{tr} ST^{-1} = \pm 2 \Leftrightarrow \lambda = 2, \ \operatorname{tr}[S, T] = \pm 2 \Leftrightarrow \lambda = \pm i$.

Jørgensen's inequality [11] tells us that $G_\lambda$ is not discrete if $|\lambda| < 1/2$.
We also note that $\lambda$ is a conjugacy invariant of the group, $G_\lambda$.

**4.1 The Riley Slice.** Riley [30] studied $\left\langle W = \begin{pmatrix} 1 & 1 \\ 0 & 1 \end{pmatrix}, K_\rho = \begin{pmatrix} 1 & 0 \\ \rho & 1 \end{pmatrix} \right\rangle$, which is conjugate to $G_\lambda$ when $\rho = 2\lambda$.

Now $\Omega(G_\lambda)/G_\lambda$ is always either a 4-punctured sphere with punctures identified in pairs or a pair of triply punctured spheres (see Section 3.3 or [27]).

Using sequences of words with real traces, Keen and Series [16] explained the fractal like boundary of

$$\mathcal{R} = \{\lambda \in \mathbb{C} \ : \ \Omega(G_\lambda)/G_\lambda \text{ is a 4-punctured sphere}\},$$

which they term the *Riley slice of Schottky space*. It is shown in Figure 6, due to D. Wright.

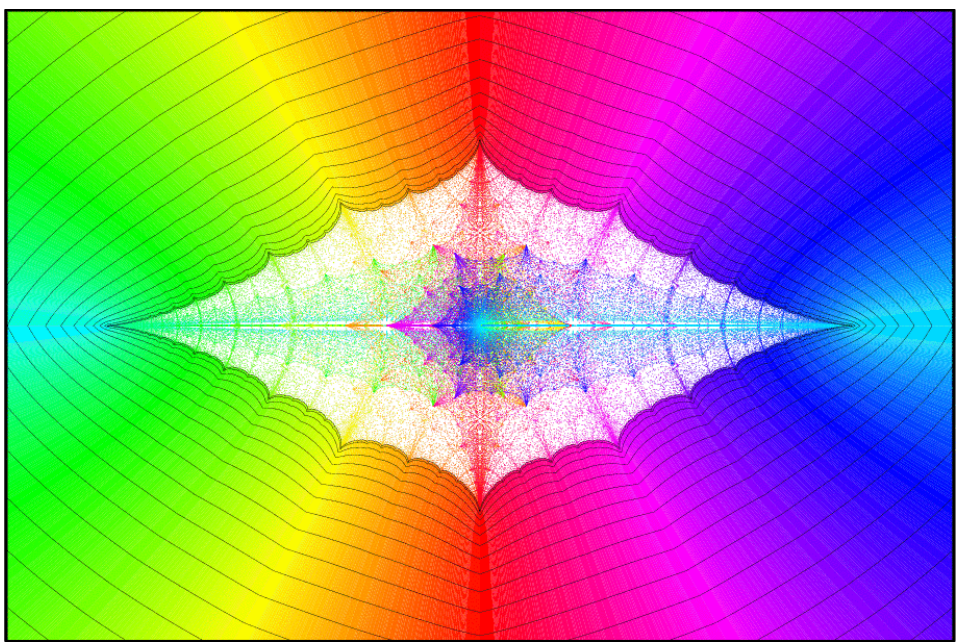

Figure 6: The Riley Slice: Courtesy of David Wright

**4.2 Spaces of groups.** We now describe various spaces of groups. For the general cases, see Marden [21], Jørgensen, Marden and Maskit [14], or Gilman [9].

For the groups generated by two parabolics that we consider, the situation is greatly simplified, since the spaces can be described as subsets of $\mathbb{C}$.

We regard $G_\lambda$ as an element of a space when $\lambda$ is.

**Definition 4.1.**

- $NSDC$-space
$$\mathcal{NSDC} = \{\lambda \in \mathbb{C} \ : \ G_\lambda \text{ is nsdc }\}$$

- *Two-parabolic $T$-Schottky* space
$$\mathcal{S} = \{\lambda \in \mathbb{C} \ : \ G_\lambda \text{ is } T\text{-Schottky }\}$$

- *Two-parabolic classical $T$-Schottky* space
$$c\mathcal{S} = \{\lambda \in \mathbb{C} \ : \ G_\lambda \text{ is Classical } T\text{-Schottky }\}$$

- *Marked two parabolic classical $T$-Schottky* space: Two parabolic generators
$$c\mathcal{S}^{\mathcal{PP}} = \left\{\lambda \in \mathbb{C} \ : \right.$$
$$\left. G_\lambda = \left\langle \begin{pmatrix} 1 & 0 \\ 1 & 1 \end{pmatrix}, \begin{pmatrix} 1 & 2\lambda \\ 0 & 1 \end{pmatrix} \right\rangle \text{ is marked Classical } T\text{-Schottky} \right\}$$

- *Marked two parabolic classical $T$-Schottky* space: Loxodromic - parabolic generators
$$c\mathcal{S}^{\mathcal{LP}} = \left\{\lambda \in \mathbb{C} \ : \right.$$
$$\left. G_\lambda = \left\langle \begin{pmatrix} 1 & 2\lambda \\ -1 & 1-2\lambda \end{pmatrix}, \begin{pmatrix} 1 & 2\lambda \\ 0 & 1 \end{pmatrix} \right\rangle \text{ is marked Classical } T\text{-Schottky} \right\}$$

- *Non-classical two-parabolic $T$-Schottky space*
$$\mathcal{NCS} = \{\lambda \in \mathbb{C} \ : \ G_\lambda \text{ is non-Classical } T\text{-Schottky }\}$$

Hence, $\mathcal{S} = c\mathcal{S} \cup \mathcal{NCS}$.

**4.3 Summary of main results** We prove (Theorem 6.7) that
$$c\mathcal{S} = c\mathcal{S}^{\mathcal{PP}}$$
and in Sections 5 and 6 show that
$$\mathcal{NSDC} \subset c\mathcal{S}^{\mathcal{LP}} \subset c\mathcal{S} \subset \mathcal{S} \subset \mathbb{C}.$$

The regions are shown in Figure 7.

Our main results, Corollary 5.13 and Theorem 6.18, can now be summarized.

**Theorem 4.2.** *Assume that $G$ is generated by two parabolic elements and is parameterized by the complex number $\lambda = x + iy$. Then*

(1) $\partial(\mathcal{NSDC}) = \{(x,y) \ : \ y^2 = 16 - 8|x|\}$;

(2) $\partial(c\mathcal{S}) = \{(x,y) \ : \ |y| = 1 - x^2/4\}$.

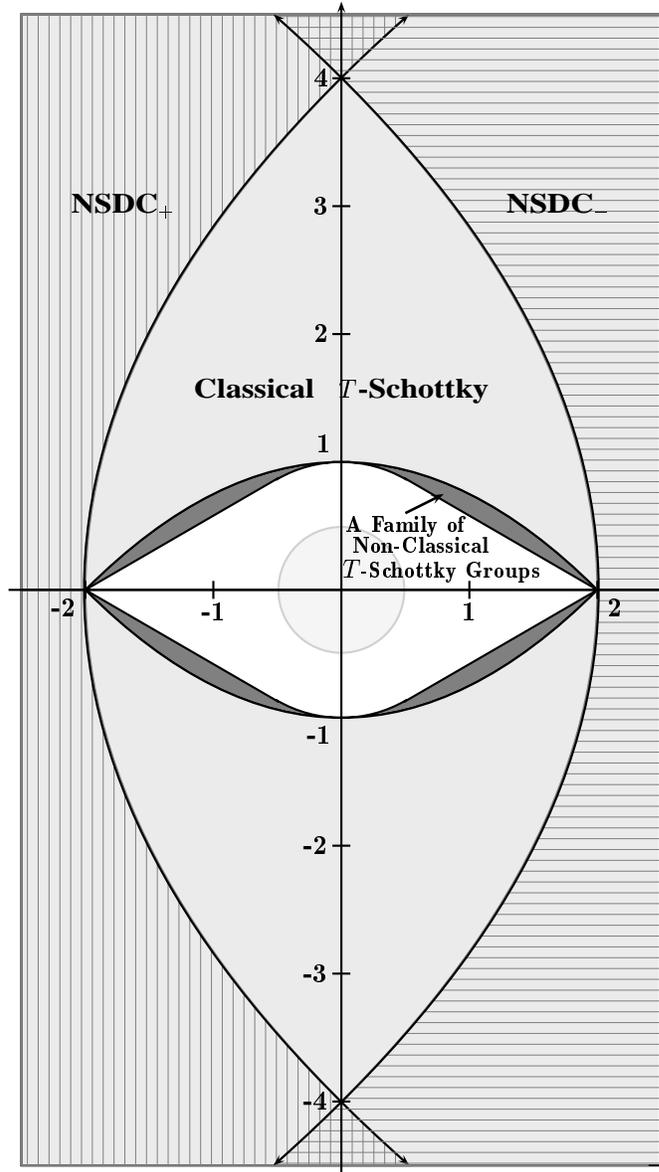

Figure 7: **Superimposed Boundary Parabolas.** Each point $\lambda \in \widehat{\mathbb{C}}$ corresponds to a two-generator group. The darkest region shows a one parameter family of non-classical $T$-Schottky groups. The line-shaded subset of the classical $T$-Schottky groups comprises the non-separating disjoint circle groups (NSDC groups). The unshaded region consists of additional non-classical $T$-Schottky groups together with degenerate groups, isolated discrete groups and non-discrete groups. Points inside the Jørgensen circle ($|\lambda| < 1/2$) are non-discrete groups.

## 5    Non-separating disjoint circle groups: NSDC Groups

Heuristically, we can think of the space of non-separating disjoint circle groups as the set of ordered six-tuples in $\widehat{\mathbb{C}}^6$ such that $(a, a', n, n', b, b')$ is the ortho-end of an nsdc-group. The boundary components of this space will be lower dimensional faces. In fact, they will be slices of nsdc space which correspond to certain limiting configurations of three circles that arise with one or more points coinciding and other extra tangencies among the three circles. Interior to the space, one can move the six points slightly and still obtain a group of nsdc type, that is, a discrete group. The boundary configurations of these slices are more rigid. If one moves a point on the boundary slightly, one sometimes obtains an nsdc group and at other times not.

In what follows, instead of working with the parameters for $\mathcal{NSDC}$, we work with the six-tuple $\sigma = (a, a', n, n', b, b')$. At the end, we write down the matrices for $A$ and $B$ in terms of the six-tuple of complex numbers. Our final conditions are given in terms of the traces of $A$, $B$ and $[A, B]$. Since the answers are given in terms of traces, these are in fact conditions on the parameters for $\mathcal{NSDC}$.

**5.1    Generators for marked nsdc groups**    We want to relate the space of two-generator nsdc groups to the space of *marked* nsdc groups. We need to sort out the effect that a change of generators has on the non-separating disjoint circles.

It is not difficult to prove the following

**Lemma 5.1.** *For a non-abelian two-generator subgroup $G$ of $PSL(2, \mathbb{C})$, the following are equivalent:*

(1) $G = \langle A, B \rangle$ *is marked nsdc;*

(2) $G = \langle B, A \rangle$ *is marked nsdc;*

(3) $G = \langle A^{-1}, B^{-1} \rangle$ *is marked nsdc;*

(4) $G = \langle A^{-1}, AB \rangle$ *is marked nsdc;*

(5) $G = \langle A^{-1}, BA \rangle$ *is marked nsdc;*

(6) $G = \langle A, A^{-1}B^{-1} \rangle$ *is marked nsdc;*

(7) $G = \langle A, B^{-1}A^{-1} \rangle$ *is marked nsdc.*

This implies

**Lemma 5.2.** *For a non-abelian two-generator subgroup $G$ of $PSL(2, \mathbb{C})$, the following are equivalent:*

(1) $G = \langle A, AB \rangle$ is marked nsdc;

(2) $G = \langle A, BA \rangle$ is marked nsdc;

(3) $G = \langle A^{-1}, A^{-1}B^{-1} \rangle$ is marked nsdc;

(4) $G = \langle A^{-1}, B^{-1}A^{-1} \rangle$ is marked nsdc.

**Lemma 5.3.** *If $G = \langle A, AB \rangle$ is marked nsdc with $A$ parabolic, then $B$ is loxodromic.*

**Proof.** Consider $A = H_A \cdot H_N$ and $AB = H_N \cdot H_{AB}$. If $C_A$, $C_N$ and $C_{AB}$ have the nsdc property, then $B = H_N H_A H_N \cdot H_{AB}$. Hence $C_{AB}$ and $B(C_{AB})$ are strictly disjoint, and $B$ maps the exterior of the first circle to the interior of the second; thus $B$ is loxodromic.

As a consequence of the above, we have

**Theorem 5.4.** *If $S$ and $T$ are parabolic, then $G = \langle S, T \rangle$ is a non-separating disjoint circle group if and only if either $\langle S, T \rangle$ or $\langle S, T^{-1} \rangle$ is marked nsdc.*

**5.2 Geometry and algebra** By noting that the half-turn $H_{[\alpha,\beta]}$ is of order two fixing $\alpha$ and $\beta$, we have

**Lemma 5.5.** *If $\alpha \neq \beta$, $\alpha \neq \infty$, and $\beta \neq \infty$, then*

$$H_{[\alpha,\beta]} = \frac{i}{\alpha - \beta} \begin{pmatrix} (\alpha + \beta) & -2\alpha\beta \\ 2 & -(\alpha + \beta) \end{pmatrix}$$

Our goal continues to be to relate the algebra and the geometry of the configuration. Since one can obtain the entries of the matrices from the six-tuple and vice-versa, conditions on the six-tuple translate into conditions on the matrices. In particular, the imposition of certain geometric configurations on the six points translate into discreteness conditions.

Corresponding to a marked group with two parabolic generators, we analyze the boundary of the slice of NSDC space when two pairs of points coalesce. In the six-tuple $\sigma = (a, a', b, b', d, d')$, we assume that $a' = b$ and $b' = d$ and term it *four-point*. Clearly, if $\sigma$ is four-point NSDC, then $C_B$ is tangent to $C_A$ at $a' = b$ and to $C_D$ at $b' = d$. Such a six-tuple together with the circles $C_A$, $C_B$ and $C_D$ is termed a *four point configuration*. Our aim is to give criteria as to when a four point configuration is indeed nsdc.

We call a point $d$ NSDC if the corresponding configuration has the nsdc property. A four point configuration with maximal tangencies between circles is a potential boundary configuration. In such a configuration, each circle is tangent to the other two. However, we shall see that not all configurations with maximal tangencies correspond to boundary points: some configurations can be pulled apart.

**Definition 5.6.** A four point configuration corresponding to a six-tuple $\sigma$ is called **extreme** if the configuration is nsdc but in any deleted neighborhood of $\sigma$, there exists a four-point $\sigma_0$ which is not nsdc and $\sigma_0'$ which is nsdc.

**5.3 Pulling Circles Apart**  We want to understand the family of circles passing through two given points in $\hat{\mathbb{C}}$. The isometric circle of $H_{[x,x']}$ has center at $(x+x')/2$ and radius $|(x-x')/2|$ (as long as neither end is at infinity). Any other circle passing through $x$ and $x'$ has its center on the perpendicular bisector of $x$ and $x'$.

Thus its center is at $c_t = \frac{x+x'}{2} + it\frac{x-x'}{2}$ for some real number $t$. The center has been pulled back $t$ units from the midpoint of the segment connecting $x$ and $x'$ along its perpendicular bisector. We let $\theta_t$ denote the angle that the radius connecting $c_t$ to $x$ or $x'$ makes with the Euclidean segment connecting $x$ and $x'$ and call $\theta_t$ the *pull-back angle* of the circle (see Figure 9).

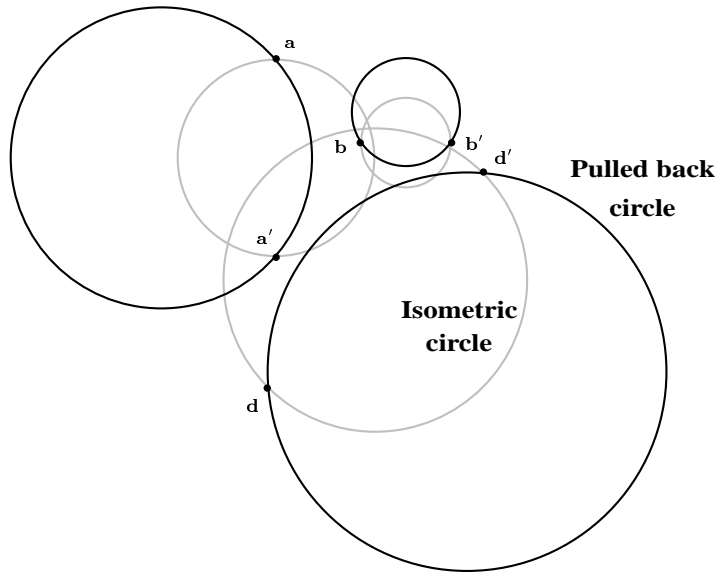

Figure 8: Isometric and pulled back circles

If we have two pairs of points, it is not difficult to analyze fully the family of pairs of circles through them; virtually all configurations of circle pairs are possible. In general, the situation is more complicated.

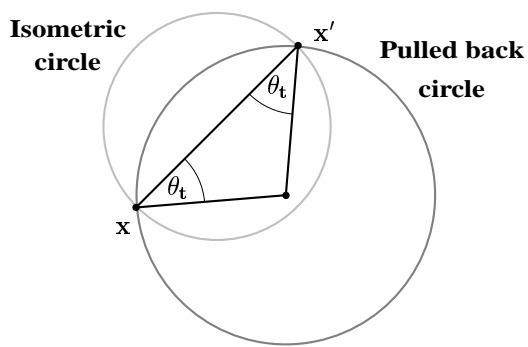

Figure 9: Pull back angle $\theta$

**5.4 Tangent Configurations** Conjugate so that $a = -2$, $a' = b = 0$ and $b' = d' = 2$. This is permissible as long as our final answers involve conjugacy invariants such as trace. The first step is to find three tangent circles, and the second step is to see which configurations of tangent circles are extreme. We want to find circles $C_A$, $C_B$ and $C_D$ with $C_A$ and $C_B$ tangent at $(0,0)$; $C_B$ and $C_D$ tangent at $(2,0)$; and $C_D$ and $C_A$ tangent at some point, call it $T$. (See Figure 10.)

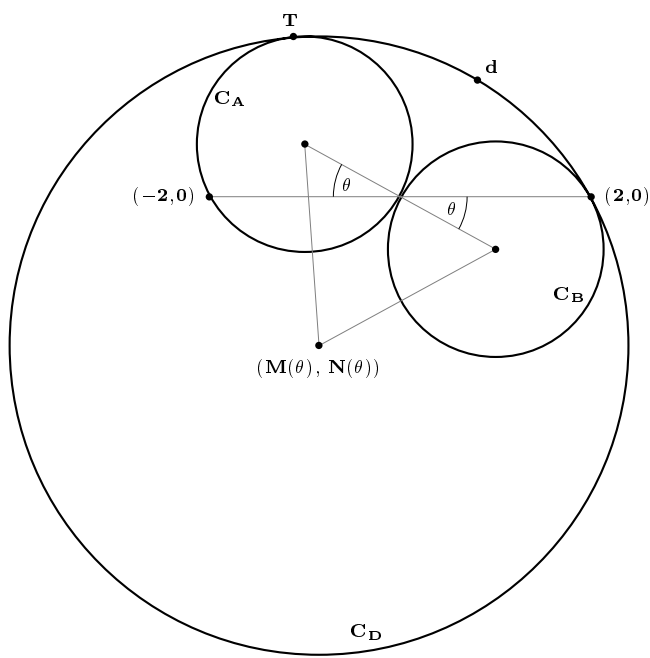

Figure 10: $\theta$ - circle configuration

The first observation is that if $C_A$ and $C_B$ are tangent at $(0,0)$, then their pull back angles must be the same but with opposite signs; thus, for this normalization, their radii are equal.

**Proposition 5.7.** (1) *Let $C_A$ and $C_B$ be pull-back circles with common pull-back angle $\theta$. For each angle $\theta$, $-\pi/4 < \theta < \pi/4$, there is a unique D-circle called the $\theta$-circle which is tangent to $C_A$ at $(2,0)$ and tangent to $C_B$.*

(2) *When $\theta \neq \pm\pi/4$ and $\theta \neq 0$, parametrically the $\theta$-circle is given by*
$x(\theta, t) = M(\theta) + R(\theta) \cdot \cos t$ *and* $y(\theta, t) = N(\theta) + R(\theta) \sin t$:

$$x(\theta, t) = \frac{2\tan\theta}{\tan\theta - (\tan\theta)^{-1}} + \frac{2\sqrt{1+\tan^{-2}\theta}}{\tan\theta - (\tan\theta)^{-1}} \cdot \cos t,$$

$$y(\theta, t) = \frac{2}{\tan\theta - (\tan\theta)^{-1}} + \frac{2\sqrt{1+\tan^{-2}\theta}}{\tan\theta - (\tan\theta^{-1})} \cdot \sin t.$$

(3) *When $\theta = 0$, the circle is:* $(x(0,t), y(0,t)) = (2\cos t, 2\sin t)$ *and when $\theta = \pm\pi/4$, the circle is:* $y = \mp(x-2)$.

**Proof.** Let the six points be $A = (-2,0)$, $A' = B = (0,0)$, $B' = D' = (2,0)$ and $D = (x,y)$. Let $L$ be the line through $(2,0)$ with slope $\tan\theta$ and $L_D$ be the line through $D$ and $(0,0)$. The cases for $\theta = 0$ and $\theta = \pm\pi/4$ are clear.

We assume the existence of a $\theta$-circle for $D$ when $C_A$ and $C_B$ both have pull-back angles $\theta$ and find formulas for its center and radius which turn out to be independent of $D$. We continue to refer to Figure 10, where $D = d$. Let $P$ be the perpendicular bisector of the segment connecting $D$ and $(2,0)$. Replacing $D$ by another point $D'$ on $L_D$ as necessary, we may assume that $L$ is not parallel to $P$. Let $(M, N)$ be the point of intersection of $L$ and $P$. Then $(M, N)$ is the center of a circle passing through $(2,0)$ and $D$, and the circle has the prescribed tangencies. We call this circle the $\theta$-circle. Note that the pull-back angle $-\theta$ is thus the same as the pull-back angle $\theta + \pi/2$.

We show that $M$ and $N$ are functions of $\theta$ and write the coordinates of the center as $(M, N) = (M(\theta), N(\theta))$. Let $c_A$ denote the center of the circle $C_A$ and $c_B$ that of $C_B$. The slope of the line passing through $(M, N)$ and $(2,0)$ is

(1) $N/(M-2) = \tan\theta$.

(2) $c_A = (-1, \tan\theta)$.

(3) $c_B = (1, -\tan\theta)$.

(4) The slope of the line segment connecting $D'$ and $D = (x, y)$ is $y/(x - 2)$.

For $x \neq 2, y \neq 0$, its perpendicular bisector passes through $(M, N)$ and $((x+2)/2, y/2)$ and has slope $-(x-2)/y$. This yields

(5) $y - 2N = -\frac{x-2}{y}(x + 2 - 2M)$.

(6) $\text{Dist}[(M, N), c_B] + R_B = R_D = \text{Dist}[(M, N), c_A] + R_A$, where $\text{Dist}[\ ,\ ]$ denotes the distance and $R_A, R_B$ and $R_D$ are the radii of the circles, $C_A, C_B, C_D$.

(7) Since $R_A = R_B = \sec\theta$, $\text{Dist}[(M, N), c_A] = \text{Dist}[(M, N), c_B]$.

(8) $(M + 1)^2 + (N - \tan\theta)^2 = (M - 1)^2 + (N + \tan\theta)^2$, whence

(9) $M = N\tan\theta$.

Conclude from $M/N = \tan\theta$ and $N/(M - 2) = \tan\theta$ that

(10) $M = M(\theta) = \frac{-2\tan^2\theta}{1-\tan^2\theta} = -\tan 2\theta \cdot \tan\theta$, and

(11) $N = N(\theta) = -\tan 2\theta$.

Writing $R = R_D$, $R = R(\theta) = \sqrt{((M-2)^2 + N^2)}$, we have

(12) $R(\theta) = |N/\sin\theta|$,

so that a point on the $\theta$-circle has coordinates given by

(13) $(x, y) = (x(\theta, t), y(\theta, t)) = (x, y) = (M + R\cos t, N + R\sin t)$ $\square$

## 5.5 Extreme four point configurations

**Definition 5.8.** A configuration of three tangent circles in $\widehat{\mathbb{C}}$ is said to be **extreme** at $d$ if the three circles $C_A, C_B$ and $C_D$ are all tangent (so that $d$ is an NSDC point) *and* if every neighborhood of $d$ contains points which are $NSDC$ and points which are not $NSDC$.

Such extreme configurations cannot be pulled apart. We find potentially extreme configurations by looking at configurations with maximal tangencies and identify extreme points by looking at the Jacobian of the mapping: $(x, y) \to (\theta, t)$.

**Proposition 5.9.** *The points in $\widehat{\mathbb{C}}$ where the Jacobian is zero are the extreme boundary points; there is one such point for each $\theta$. All other points either have a totally NSDC neighborhood or a totally non-NSDC neighborhood. The Jacobian is non-zero except at the point $z_0 = (x_0, y_0)$ where*

(4) $$x_0 = 2 - 4\frac{\sin^2\theta}{1 + \sin^2 2\theta} \quad \text{and} \quad y_0 = 8\frac{\sin^3\theta \cdot \cos\theta}{1 + \sin^2 2\theta}.$$

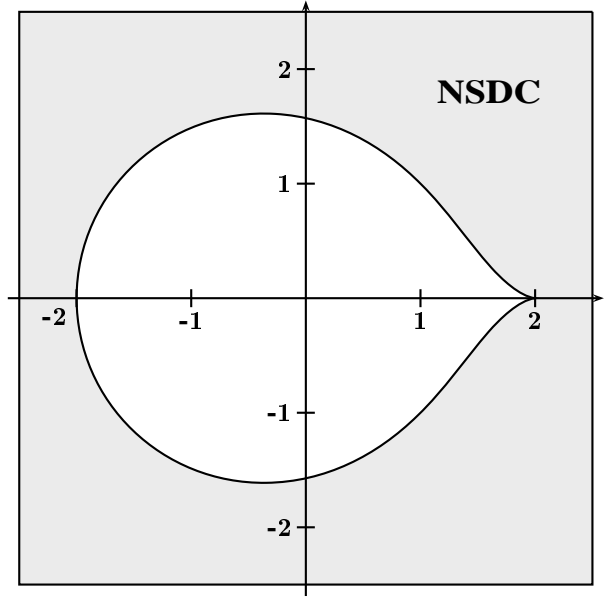

Figure 11: **The NSDC teardrop.** The teardrop represents the plot of the critical points $(x_0, y_0)$. It is the same as the envelope of the family of $\theta$-circles. In particular, all other points on a given $\theta$-circle have a totally NSDC neighborhood. Points $d$ in the complex plane that lie exterior to the teardrop are of non-separating disjoint circle type.

**Proof.** Compute partial derivatives.

(1) $x_\theta = M_\theta + R_\theta \cos t$, $x_t = -R \sin t$ and $y_\theta = N_\theta + R_\theta \sin t$, $y_t = R \cos t$.

Computing the Jacobian and setting it equal to zero yields

(2) $RR_\theta + RM_\theta \cos t + RN_\theta \sin t = 0$.

Using the facts that $RR_\theta = M_\theta(M-2) + NN_\theta$, $\cos t = (x-M)/R$ and $\sin t = (y-N)/R$, simplify the Jacobian to obtain that

(3) the Jacobian is $0$ if and only if $y/(x-2) = -(M_\theta/N_\theta)$.

(4) Compute these partials to conclude that the critical points of the Jacobian are precisely those points that satisfy: $y/(x-2) = -\sin 2\theta$.

Namely writing, $M = -\tan 2\theta \cdot \tan \theta$, we have $M_\theta = -(2\sec^2 2\theta \tan \theta + \sec^2 \theta \tan 2\theta)$. Since $N = N(\theta) = -\tan 2\theta$, $N_\theta = -2\sec^2 2\theta$. Thus

$$-\frac{M_\theta}{N_\theta} = \frac{2\sec^2 2\theta \tan \theta + \sec^2 \theta \tan 2\theta}{-2\sec^2 2\theta} = -\sin 2\theta.$$

(5) Combine the fact that the distance from $(x, y)$ to $(M, N)$ is $R$ with $y/(x-2) = -\sin 2\theta$ to obtain the coordinates for extreme points (equation 4).

(6) If $(x_0, y_0)$ are given by the formula (4), calculate that $y_0/(x_0 - 2) = -\sin 2\theta$ so that $(x_0, y_0)$ is extreme.

It is straightforward to show that for each pair $(\theta, t)$, the point $d = x(\theta, t) + iy(\theta, t)$ or the six-tuple $(-2, 0, 2, 0, 2, d)$ corresponds to an NSDC group where the three circles are tangent. For each $\theta$, the points $(-2, 0, 2, 0, 2, d)$ with $d$ exterior to the $\theta$-circle are NSDC. $\square$

**5.6 Matrix Equivalents** By Lemma 5.5,

$$H_{[-2,0]} = \begin{pmatrix} -i & 0 \\ i & i \end{pmatrix}$$

$$H_{[0,2]} = \begin{pmatrix} -i & 0 \\ -i & i \end{pmatrix}$$

$$H_{[2,d]} = \frac{i}{2-d}\begin{pmatrix} (d+2) & -4d \\ 2 & -(d+2) \end{pmatrix}$$

and hence

$$\gamma_0 = H_{[-2,0]}H_{[0,2]}H_{[2,d]} = \frac{i}{2-d}\begin{pmatrix} d+2 & -4d \\ -2d-2 & 7d-2 \end{pmatrix}$$

$$A = H_{[-2,0]}H_{[0,2]} = \begin{pmatrix} 1 & 0 \\ -2 & 1 \end{pmatrix}$$

$$B = H_{[0,2]}H_{[2,d]} = \frac{1}{2-d}\begin{pmatrix} d+2 & -4d \\ d & -3d+2 \end{pmatrix}$$

and

$$tr(A) = tr(B) = 2, \ tr(\gamma_0) = 8di/(2-d) \text{ and } tr(AB) = 2 + 8d/(2-d).$$

As before, $\operatorname{tr}(X)$ denotes the trace of the matrix $X$ and $[X, Y]$ denotes the multiplicative commutator of two matrices so, that $\gamma_0^2 = [A, B]$. We have shown

**Proposition 5.10.** *If $tr(A) = 2$, $tr(B) = 2$, and $tr(AB) - 2 = 2\frac{4d}{2-d}$, where $d$ lies exterior to the NSDC teardrop, then $G = \langle A, B \rangle$ is discrete.*

*If $tr(A) = 2$, $tr(B) = 2$, and $tr[A, B] = 8di/(2-d)$, where $d$ lies exterior to the NSDC teardrop, then $G = \langle A, B \rangle$ is discrete.*

This motivates the following re-normalization: The map $z \to 4z/(z-2)$ sends the six-tuple $(-2, 0, 0, 2, 2, d)$ to the six-tuple $(2, 0, , 0, \infty, \infty, \lambda)$, where $\lambda = 4d/(d-2)$.

When the parameter $d$ is replaced by $\lambda$, the teardrop is transformed into a parabola.

**5.7 Extreme four point configurations via** $(2, 0, 0, \infty, \infty, \lambda)$   We note that

$$S = H_{[2,0]}H_{[0,\infty]} = \begin{pmatrix} 1 & 0 \\ 1 & 1 \end{pmatrix} \quad \text{and} \quad T = H_{[0,\infty]}H_{[\infty,\lambda]} = \begin{pmatrix} 1 & -2\lambda \\ 0 & 1 \end{pmatrix}$$

and hence
$$G = \langle S, T \rangle = G_{-\lambda}.$$

Motivated by the above, we term the marked group $G_\lambda$ of NSDC$_-$ type if the six-tuple $(2, 0, 0, \infty, \infty, \lambda)$ is NSDC and of NSDC$_+$ type if the six-tuple $(2, 0, 0, \infty, \infty, -\lambda)$ is NSDC. By Theorem 5.4, $G_\lambda$ is NSDC if and only if the marked group is either NSDC$_-$ or NSDC$_+$.

We write $\lambda = |\lambda|e^{i\omega}$, as in Section 4.

The four point configuration for $(2, 0, 0, \infty, \infty, \lambda)$ involves a circle $C_1$ through $0$ and $2$ together with a line $L_2$ tangent to $C_1$ at $0$ and a line $L_3$ through $\lambda$. $C_1$ may be parameterized by specifying the angle $\phi$ between the $x$-axis and the radial vector through $0$, measured in the anti-clockwise direction. Thus $-(\pi/2) < \phi < (\pi/2)$. The radius of $C_1$ is then $R(\phi) = 1/\cos\phi$. Once $C_1$ is specified, $L_2$ and $L_3$ are determined. We let $D_L(\phi)$ be the signed distance between $L_2$ and $L_3$, the positive direction being the radial vector from $0$ to the center of $C_1$.

The condition that the configuration be an nsdc configuration for this normalization is that $L_2$ and $L_3$ be sufficiently far apart to contain $C_1$ and that $C_1$ actually lie between these two lines.

The following is then the criterion we are seeking

**Lemma 5.11.** $C_1$, $L_2$ and $L_3$ form an NSDC configuration for $(2, 0, 0, \infty, \infty, \lambda)$

$$\Leftrightarrow \quad D_L(\phi) \geq 2R(\phi)$$
$$\Leftrightarrow \quad |\lambda|\cos(\omega - \phi)\cos(\phi) \geq 2$$

**Proof.** The condition $|D_L(\phi)| \geq 2R(\phi)$ guarantees that $L_2$ and $L_3$ are sufficiently far apart to contain $C_1$, and $D_L(\phi) > 0$ ensures that $C_1$ lies between $L_2$ and $L_3$. The equivalent condition follows on noting that

$$D_L(\phi) = \lambda \cdot e^{i\phi} \quad \text{and} \quad R(\phi) = \frac{1}{\cos\phi}. \qquad \square$$

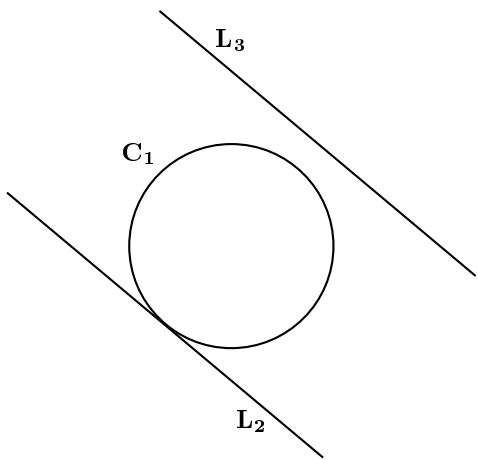
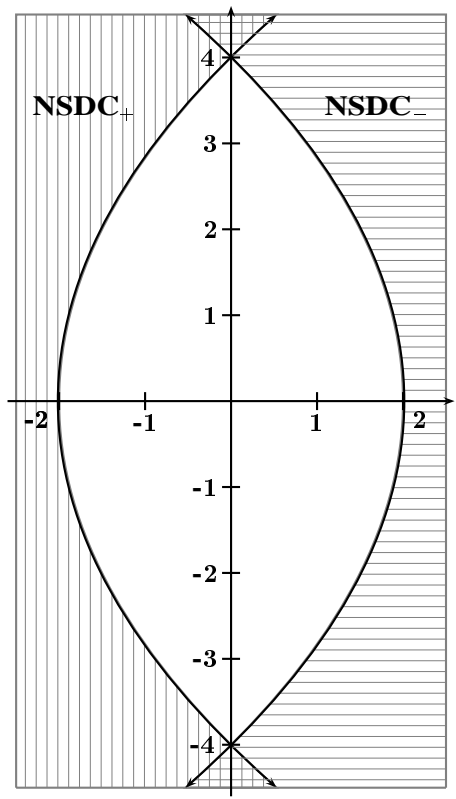

Figure 12: Normalization through $\infty$

Figure 13: Boundary Parabolas - NSDC

Finally, we can identify when $G_\lambda$ is nsdc.

**Theorem 5.12.**

$$G_\lambda \text{ is NSDC}_- \quad \Leftrightarrow \quad |\lambda|[1 + \cos\omega] \geq 4$$
$$G_\lambda \text{ is NSDC}_+ \quad \Leftrightarrow \quad |\lambda|[1 - \cos\omega] \geq 4.$$

**Proof.** By Lemma 5.11, $G_\lambda$ is NSDC$_-$ $\Leftrightarrow$

$$|\lambda|[\cos\omega + \cos(\omega - 2\phi)] \geq 4.$$

Maximizing the left hand side as a function of $\phi$ gives the desired result; the maximum occurs when $\phi = \omega/2$. $\square$

Points on the boundary of $\mathcal{NSDC}$-space correspond to groups with maximal tangencies that are extreme. We have shown

**Corollary 5.13.**

$$G_\lambda \text{ is NSDC} \quad \Leftrightarrow \quad |\lambda|[1 + |\cos\omega|] \geq 4.$$

*The boundary curve $|\lambda|[1 + |\cos\omega|] = 4$ gives the boundary of NSDC space in the polar coordinates $(|\lambda|, \omega)$, and this boundary curve is the piecewise-parabola whose Euclidean coordinates $(x,y)$ satisfy $y^2 = 16 - 8|x|$, as shown in Figure 13.*

Finally, we put Theorem 5.12 into a conjugacy invariant form.

**Theorem 5.14.** *Assume that $G = \langle S, T \rangle$ with $tr(S) = tr(T) = 2$. Then $G$ is NSDC if and only if*

$$|tr(ST) - 2| + \big|Re[tr(ST) - 2]\big| \geq 8.$$

**Proof.** The result follows by recalling that $\lambda = (1/2)[tr(ST) - 2]$. □

Since an nsdc-group is discrete, we have the following discreteness criteria.

**Corollary 5.15.** *Assume that $G = \langle S, T \rangle$ with $tr(S) = tr(T) = 2$. Then*

$$\big|tr(ST) - 2\big| + |Re[tr(ST) - 2]| \geq 8 \implies G \text{ is discrete.}$$

## 6   $T$-Schottky Groups

In what follows, we continue the goal of relating the geometry of a configuration to its algebra. One can view our results as describing precisely when requiring tangencies forces a group to be a boundary group of $T$-Schottky space and when tangencies can be *pulled apart* so that the group is an interior group of $T$-Schottky space.

Section 6.2 develops the main geometric insights.

Preliminary to the geometric analysis we describe the effects of marking generators.

**6.1   Change of generators**   In this subsection, we prove a series of lemmas, the first of which is an easy observation. These lead up to the proof of Theorem 6.7 about the relation between classical $T$-Schottky groups and marked classical $T$-Schottky groups for groups generated by two-parabolics.

**Lemma 6.1.** *For a subgroup $G$ of $PSL(2, \mathbb{C})$, the following are equivalent:*

(1) $G = \langle A, B \rangle$ *is marked classical $T$-Schottky;*

(2) $G = \langle B, A \rangle$ *is marked classical $T$-Schottky;*

(3) $G = \langle A, B^{-1} \rangle$ is marked classical $T$-Schottky.

We say that a fixed point of a parabolic transformation is represented on the boundary of a $T$-Schottky domain if there is an element of the group that maps the parabolic fixed point to a point on the boundary of the domain. We prove the following refinement of a well-known result.

**Lemma 6.2.** *If $G$ is a classical $T$-Schottky group containing a parabolic $T$ and $F$ is a classical $T$-Schottky domain for $G$, then the fixed point of $T$ is represented on the boundary of $F$ and two sides of $F$ are tangent at that point.*

**Proof.** Normalize so that $T$ fixes infinity and assume that infinity is not represented on the boundary of $F$. Since $\infty$ is a fixed point of $T$, $\infty$ is not an interior point of $F$. There is then a boundary circle $C$ separating infinity from the other bounding circles. Let $C'$ be the bounding circle that is paired to this $C$ and $g_C$ the pairing. Specifically, choose $g_C \in G$ to send the exterior of $C'$ to the interior of $C$. We consider the image of $F$ paired to this circle, that is, $g_C(F)$. If $\infty$ is on the boundary of $g_C(F)$, we are done. Otherwise, $g_C(F)$ is again a $T$-Schottky domain for $G$ with $\infty$ in its exterior. Repeat this process. Since $\infty$ is a limit point, these nested images of $F$ accumulate at $\infty$. However, $G$ contains the translation $T$, so at some stage the image of $F$ will overlap itself.

For the second part, (re)-normalize so that infinity is on the boundary and note that if only one bounding circle passes through infinity, then there would be an element $g \in G$ which is not the identity with $g(F) \cap F \neq \emptyset$. □

If $F$ is a classical $T$-Schottky domain for a group $G$, then for any $g \in G$, either $g(F) \cap F = \emptyset$ or $g(F) = F$. We say that $F$ and $g(F)$ are tangent if there are circles $C$ and $D$ bounding $F$ with $C$ and $g(D)$ tangent.

**Lemma 6.3.** *Let $F$ be a classical $T$-Schottky domain and $S$ a side pairing transformation. If $F$ and $S^N(F)$ are tangent at $\eta$ for some $N$ with $|N| \geq 2$, then $S$ fixes $\eta$.*

**Proof.** Since $F$ is a classical $T$-Schottky domain, there are circles $C$ and $C'$ on $\partial(F)$ with $S(C) = C'$, where $S$ maps the exterior of $C$ to the interior of $C'$. If $\eta$ lies interior to $\overline{F}$, then $S(\eta)$ is interior to $C'$, but $F$ and $S^N(F)$ are disjoint. Thus we must have $\eta \in C \cap C'$. Since $|N| \geq 2$, it follows that all $S^n(\overline{F})$ are tangent at $\eta$ for all integers $n$. The $S$ orbit of a circle through $\eta$ accumulates at a fixed point of $S$, so $\eta$ must be such a fixed point. □

**Lemma 6.4.** *Suppose $G = \langle S, T \rangle$ is marked classical $T$-Schottky. Then $S^N T$ and $(ST)^N T$ are loxodromic for $|N| \geq 2$.*

**Proof.** Let $F$ be a classical $T$-Schottky domain for the marked group $G$ with $T$ pairing sides $C$ and $C'$. Then $S^N T$ maps the exterior of $C$ to the interior of $S^N T(C)$, so the region bounded by these two circles is itself a $T$-Schottky domain for the cyclic group generated by $S^N T$. Therefore, if $S^N T$ is parabolic, the two circles $C$ and $S^N T(C)$ must be tangent at a point $\eta$. Thus $\overline{F}$ and $S^N(\overline{F})$ are tangent at $\eta$. But then by Lemma 6.3, $\eta$ is fixed by $S$. Therefore, both $C$ and the sides of $F$ paired by $S$ are tangent at $\eta$. However, since there are no elliptic elements in $T$-Schottky groups, such a configuration is impossible. Therefore, $S^N T$ must be loxodromic.

Likewise, if $(ST)^N T$ is parabolic, then $C$ and $(ST)^N T(C)$ must be tangent at a point $\zeta$. Applying Lemma 6.3 to the domain bounded by $C$ and $ST(C)$, which is a classical $T$-Schottky domain for the cyclic group $\langle ST \rangle$, we see that $ST$ fixes $\zeta$. Hence $C$, $T(C)$ and $ST(C)$ are tangent at $\zeta$, which cannot happen. Thus $(ST)^N T$ must be loxodromic. □

In contrast to the purely loxodromic case, where Chuckrow [5] has shown that every Schottky group is a marked Schottky group on any set of free generators, we have

**Theorem 6.5.** *Suppose that $G = \langle A, B \rangle$ is a marked classical $T$-Schottky group. If $G$ can also be generated by parabolic elements $S$ and $T$ then, up to conjugacy, interchange of generators and replacing a generator by its inverse, we have*

*either*    (i) $S = A$ *and* $T = B$    *or*    (ii) $S = A$ *and* $T = AB$.

**Proof.** Let $F$ be a classical $T$-Schottky domain for the marked group $G$. Observe that $F$ can have at most six points of tangency; that is, the four circles that bound $F$ can have a total of at most six points of tangency. If any two circles are tangent at a point, none of the other four circles can be tangent at that point, and each circle can be tangent to at most two others.

Since the fixed points of both $S$ and $T$ are represented on the boundary of $F$ as points of tangency by Lemma 6.2 (see Figure 17), there are words $W_1(A, B)$ and $W_2(A, B)$ with fixed points on the boundary of $F$ that are conjugate to $S$ and $T$ respectively. However, there may not necessarily be a single element of the group that conjugates $W_1(A, B)$ to $S$ and $W_2(A, B)$ to $T$.

Interchanging $S$ and $T$ if necessary, we may conjugate so that the fixed point of $S$ lies on the boundary of $F$ and is fixed by the shortest possible word in $A$ and $B$.

Elements of the group send points of tangency of $F$ to points of tangency of the image of $F$ and parabolic fixed points to parabolic fixed points. We observe that the generators $A$ and $B$ must each permute the parabolic fixed points on the boundary of $F$: if $A$ identifies $L$ and $M$ and $B$ identifies $C$ and $D$, we assume $p$ is a parabolic fixed lying on $L$ and consider the possibilities for its images. Since $A(L) = M$ if $A(p)$ lies on $M$, $L$ is either tangent to $C$,$D$ or $M$ at $p$; and by Lemma 6.2, $M$ is tangent at $A(p)$ to one of $L$, $C$ or $D$. Considering all possible cases of tangency points under the action of the group, we conclude that it suffices to consider the three cases in which the length of $S$ is at most three. Now if $U$ and $V$ are associated primitive elements of the free group on two generators, then the only primitive elements associated to $U$ are of the form $U^\alpha V^\epsilon U^\beta$, where $\alpha, \beta = 0, \pm 1, \pm 2, \ldots$ and $\epsilon = \pm 1$ (p.169 of [18]). We apply this to the cases where the length of $S$ is at most three.

Length 1. $S$ is a word of length one in $A$ and $B$. We may assume that $S = A$ and note that, up to conjugacy by a power of $A$, we get $T = A^n B^{\pm 1}$. By Lemma 6.4, either $n = 0$ or $n = \pm 1$. The theorem follows after a re-normalization.

Length 2. $S$ is a word of length two in $A$ and $B$. We may assume that $S = AB$ and $T = (AB)^n B$. By Lemma 6.4, either $n = 0$ or $n = \pm 1$. If $n = 0$ or $n = -1$, then up to conjugacy, $T$ has length one. If $n = 1$, then Lemma 6.4 shows that $T$ is loxodromic. Both cases are contrary to assumption.

Length 3. $S$ is a word of length three in $A$ and $B$. We may assume that $S = A^2 B$. However, by Lemma 6.4, $S$ is then loxodromic, contrary to assumption. $\square$

**Lemma 6.6.** *If $\langle T, B \rangle$ is marked classical $T$-Schottky with $T$ and $TB$ parabolic, then $\langle T, TB \rangle$ is marked classical $T$-Schottky.*

**Proof.** Normalize so that $T(\infty) = \infty$ and $TB(0) = 0$ and hence $B(0) = T^{-1}(0)$. A $T$-Schottky domain $F$ for $\langle T, B \rangle$ is then bounded by $L, M, C$ and $D$, where $T$ pairs sides $L$ and $M$ that meet at infinity and $B$ pairs circles $C$ and $D$ inside the domain bounded by $L$ and $M$. By Lemma 6.2 and Theorem 6.5, $0$ lies on the boundary of $F$.

There are the following three cases.

(1) $0$ on $L$, $T(0)$ on $M$.

 Since the Euclidean line $L$ separates $T^{-1}(L)$ and $M$, we must have $B(0) \neq T^{-1}(0)$.

(2) $0$ is interior to the domain bounded by $L$ and $M$ with $C$ and $D$ tangent at $0$.

 In this case, $T^{-1}(0)$ is interior to the region bounded by $T^{-1}(L)$ and $L$, and $B(0)$ is on $D$. Again, $B(0) \neq T^{-1}(0)$.

(3) $C$ is tangent to $M$ at $0$, and $D$ is tangent to $L$ at $T^{-1}(0)$.

We have $B(0) = T^{-1}(0)$.

In case (3), we let $L_1$ be a line orthogonal to the line joining the centers of $C$ and $D$ and separating the interiors $C$ and $D$. Such a line exists, since $C$ and $D$ are at worst tangent. If $M_1 = T(L_1)$ and $D_1 = T(D)$, then $T$ pairs $L_1$ and $M_1$ and $TB$ pairs $C$ and $D_1$. Thus $\langle T, TB \rangle$ is marked classical $T$-Schottky. □

Theorem 6.7 is a direct consequence of Lemmas 6.1 and 6.6, together with Theorem 6.5.

**Theorem 6.7.** *If $S$ and $T$ are parabolic, then $G = \langle S, T \rangle$ is marked classical $T$-Schottky if and only if $G$ is classical $T$-Schottky.*

**6.2 Geometry of marked two parabolic generator $T$-Schottky groups.** In what follows, we consider a marked group $G$ with two parabolic generators, $S$ and $T$, and parameter (Section 4) the complex number $\lambda$. We write $G = G_\lambda$ with $\lambda = |\lambda|e^{i\omega}$. Replacing $T$ by $T^{-1}$, if necessary, we may assume by Lemma 6.1 that $0 \leq \omega < \pi$.

The parameter $\lambda = |\lambda|e^{i\omega}$ that comes from the normalized generators determines certain geometric parameters which we describe below. Our goal is to determine how the geometric parameters determine $|\lambda|$ and $\omega$, the algebraic parameters defined in Section 4.

We consider $S$ pairing circles $C_1$ and $C_2$ and mapping the exterior of $C_1$ to the interior of $C_2$ and the translation $T$ pairing distinct lines $L_1$ to $L_2$. Since $S$ and $T$ are parabolic, $C_1$ is tangent to $C_2$ at $0$, the fixed point of $S$, just as $L_1$ is tangent to $L_2$ at infinity, the fixed point of $T$.

The configuration $C_1$, $C_2$, $L_1$ and $L_2$ determines a classical $T$-Schottky domain provided that $C_1$ and $C_2$ lie in the region bounded by $L_1$ and $L_2$. If they do, then $G = \langle S, T \rangle$ is discrete and a classical $T$-Schottky group. We thus assume that $0$ lies in the region bounded by $L_1$ and $L_2$.

**Definition 6.8.** We call a classical $T$-Schottky group $G_\lambda$ and the corresponding classical $T$-Schottky configuration **extreme** if every deleted neighborhood of the parameter $\lambda$ contains points for which the group is not classical $T$-Schottky and points for which it is classical $T$-Schottky.

Note that if either $L_1$ is not tangent to either circle or $L_2$ is not tangent to either circle, then the configuration is not extreme. Indeed if $\lambda_0$ is sufficiently close to $\lambda$, one can easily find $L_1'$ paired to $L_2'$ by $T(z) = z + 2\lambda_0$ so that $C_1$, $C_2$ lie strictly within the strip bounded by $L_1'$ and $L_2'$.

The situation is more subtle if $L_1$ is tangent to $C_1$ and $L_2$ is tangent to $C_2$. Clearly, extreme configurations are of this type. However, it turns out that not all such configurations are extreme. That is, generically one expects an extreme configuration to have six tangencies; but not all six tangency configurations are extreme, and a candidate for an extreme configuration may have only four tangencies.

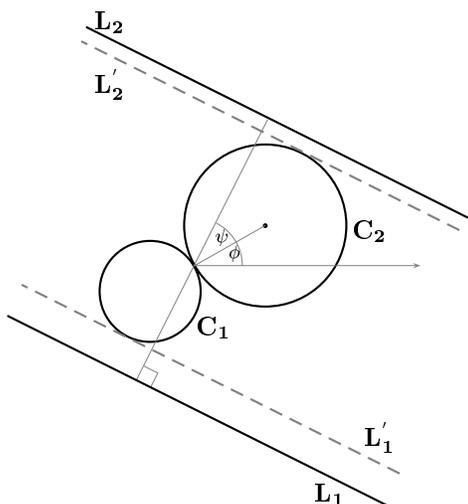
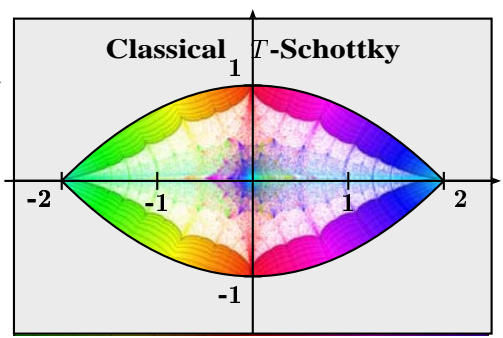

Figure 14: $T$-Schottky configuration: Two parabolic generators.

Figure 15: Boundary Parabolas - Classical $T$-Schottky

### 6.2.1 The parameters $t$, $\phi$ and $\psi$.

We want to assign *geometric parameters* to the set of four circles.

Let the two circles $C_1$ and $C_2$ be tangent at $0$ and assume that the line connecting their centers makes an angle $\phi$ with the $x$-axis. Let $\tau = te^{i\phi}$ and assume that the center of $C_1$ is $\alpha = -\tau$.

Applying Lemma 8.3 to this situation with $\nu = 1$ so that $S$ is parabolic, we obtain

**Corollary 6.9.** *If the center of $C_1$ is at $-\tau$, the center of $C_2$ is at $\tau/(\tau + \overline{\tau} - 1)$. The exterior of $C_1$ is mapped to the interior of $C_2$ if and only if $\tau + \overline{\tau} - 1 > 0$.*

Observe that by Corollary 6.9, in order for the exterior of $C_1$ to be mapped to the interior of $C_2$ it is necessary that $|\phi| \leq \pi/2$.

Assume that the perpendicular to $L_1$ and $L_2$ makes an angle $\phi + \psi$ with the $x$-axis, $|\psi| \leq \pi/2$.

**Definition 6.10.** Let $D_L = D_L(\phi, \psi)$ denote the Euclidean distance between $L_1$ and $L_2$.

**Definition 6.11.** Let $D_C = D_C(t, \phi, \psi)$ be the distance between the parallel lines $L'_1$ and $L'_2$, where $L'_1$ is tangent to $C_1$ and parallel to $L_1$ and $L'_2$ is tangent to $C_2$, with $L'_1$ and $L'_2$ chosen so that $C_1$ and $C_2$ lie between $L'_1$ and $L'_2$.

We emphasize that *the condition that we have a $T$-Schottky configuration is thus precisely that $D_C \leq D_L$.*

In what follows, we progressively describe how best to choose the *geometric parameters* $t, \phi, \psi$. Essentially, we minimize $D_C - D_L$ as a function of $t$, then $\phi$, and then $\psi$.

**Lemma 6.12.**

$$\begin{aligned} D_C &= t\left[1 + \frac{1}{(2t\cos\phi - 1)}\right][1 + \cos\psi] \\ D_L &= |2\lambda \cdot e^{i(\phi+\psi)}| \end{aligned}$$

**Proof.** The result is a straightforward computation. □

**Definition 6.13.** $t_0 = 1/\cos\phi$.

**Lemma 6.14.** *For fixed $\phi$ and $\psi$, there is a $T$-Schottky configuration if and only if $D_C(t_0, \phi, \psi) \leq D_L(\phi, \psi)$, where $t_0 = 1/\cos\phi$. When $t = t_0$, $D_C(t_0, \phi, \psi) = (2(1 + \cos\psi))/\cos\phi$, this value of $D_C$ is a minimum, and the circles have the same radius.*

**Proof.** By Lemma 6.12,

$$D_C = \left[(2t\cos\phi - 1) + \frac{1}{(2t\cos\phi - 1)} + 2\right]\frac{(1 + \cos\psi)}{2\cos\phi};$$

and the result follows. □

**Definition 6.15.**
$$\phi_0 = \left\{\begin{array}{ll} \frac{1}{2}(\omega - \psi) & \text{when} \quad |\omega - \psi| \leq \frac{\pi}{2} \\ \frac{1}{2}(\omega - \psi - \pi) & \text{when} \quad \omega - \psi \geq \frac{\pi}{2} \end{array}\right\} \quad \text{In both cases,} \quad |\phi_0| \leq \frac{\pi}{4}.$$

**Lemma 6.16.** *For fixed $\psi$, there is a $T$-Schottky configuration if and only if $D_L(\phi_0, \psi) \geq D_C(t_0, \phi_0, \psi)$, where $\phi_0$ satisfies 6.15.*

**Proof.** By Lemma 6.14, there is a $T$-Schottky configuration if and only if

$$|2\lambda \cdot e^{i(\phi+\psi)}| \geq \frac{2(1+\cos\psi)}{\cos\phi}.$$

Further,

$$|2\lambda \cdot e^{i(\phi+\psi)}| \geq \frac{2(1+\cos\psi)}{\cos\phi}$$

if and only if

$$|\lambda||\cos(\omega-\psi) + \cos(\omega-\psi-2\phi)| \geq 2(1+\cos\psi);$$

and the result follows by maximizing the left hand side as a function of $\phi$.   □

Utilizing the above, we explicitly characterize classical $T$-Schottky groups on two parabolic generators.

**Theorem 6.17.** *$G_\lambda$, $0 \leq \omega < \pi$, is marked classical $T$-Schottky if and only if*

$$|\lambda|(1+\sin\omega) \geq 2.$$

*If the above inequality holds, then a $T$-Schottky configuration may be obtained by choosing*

$$\begin{aligned}
\psi &= \pm(\pi/2), \\
\phi &= \omega/2 - \pi/4, \quad \text{and} \\
t &= |1 + i\tan(\omega/2 - \pi/4)|.
\end{aligned}$$

*Further, the two choices of triples above, $(\pi/2, \phi, t)$ and $(-\pi/2, \phi, t)$, are the only choices of parameters $(\psi, \phi, t)$ guaranteed to give a classical $T$-Schottky configuration for every $\omega$ ($0 \leq \omega < \pi$).*

*If $\omega = 0$, then every classical $T$-Schottky configuration can be obtained by choosing some $\psi$, then setting $\phi = \phi_0$ according to Definition 6.15.*

**Proof.** From the proof of Lemma 6.16, $G_\lambda$ is marked classical $T$-Schottky if and only if

$$|\lambda|[1 + |\cos(\omega-\psi)|] \geq 2(1+\cos\psi).$$

But

$$|\lambda|[1 + |\cos(\omega-\psi)|] \geq 2(1+\cos\psi)$$

$$\iff \begin{cases} |\lambda|\cos^2(\frac{\omega}{2} - \frac{\psi}{2}) \geq 2\cos^2\frac{\psi}{2} \quad \text{and} \quad -\frac{\pi}{2} \leq \omega - \psi \leq \frac{\pi}{2} \\ \quad \text{or} \\ |\lambda|\sin^2(\frac{\omega}{2} - \frac{\psi}{2}) \geq 2\cos^2\frac{\psi}{2} \quad \text{and} \quad \frac{\pi}{2} < \omega - \psi \leq \frac{3\pi}{2} \end{cases}$$

$$\iff \begin{cases} |\lambda|[\cos(\frac{\omega}{2}) + \sin(\frac{\omega}{2})\tan(\frac{\psi}{2})]^2 \geq 2 \quad \text{and} \quad -\frac{\pi}{2} \leq \omega - \psi \leq \frac{\pi}{2} \\ \quad \text{or} \\ |\lambda|[\sin(\frac{\omega}{2}) - \cos(\frac{\omega}{2})\tan(\frac{\psi}{2})]^2 \geq 2 \quad \text{and} \quad \frac{\pi}{2} < \omega - \psi \leq \frac{3\pi}{2}. \end{cases}$$

Provided $\omega \neq 0$, the maximum of the left hand side as a function of $\psi$ occurs precisely when $\psi = \pm(\pi/2)$.

If $\omega = 0$, then the left hand side is independent of $\psi$, so any $\psi$ may be chosen. □

Figure 17 illustrates the generic extreme configuration (i.e., $\omega \neq 0$).

As a consequence, we can explicitly describe the boundary of classical $T$-Schottky space with two parabolic generators simply as a portion of a parabola (Figures 15, 7). We have $x + iy = |\lambda|e^{i\omega}$.

**Theorem 6.18.** $G_\lambda$ *lies on the boundary of classical T-Schottky space* $\iff$

$$\lambda = (2e^{i\omega})/(1 + |\sin \omega|)$$

*and thus, eliminating $\omega$,* $\iff$

$$\lambda = x + iy \text{ with } |y| = 1 - x^2/4$$

.

Finally, we put Theorem 6.17 into a conjugacy invariant form,

**Theorem 6.19.** *If $G = \langle S, T \rangle$ with $tr(S) = tr(T) = 2$, then $G$ is classical T-Schottky if and only if*

$$|tr(ST) - 2| + |Im[tr(ST)]| \geq 4.$$

**Proof.** The result follows by recalling that $\lambda = (1/2)[tr(ST) - 2]$. □

This translates to the following sufficient condition for discreteness.

**Corollary 6.20.** *If $G = \langle S, T \rangle$ with $tr(S) = tr(T) = 2$,*

$$|tr(ST) - 2| + |Im[tr(ST)]| \geq 4$$

*implies that $G$ is discrete.*

**6.3 Generically not pinching** A point in the two-parabolic $T$-Schottky space represents a surface that has been pinched precisely when a new conjugacy class of primary parabolic elements appears. That is, boundary points come from additional pinching only when the pinched groups have a third conjugacy of primary parabolics. This third conjugacy class may or may not be primitive.

A corollary of our results is that points on the boundary $\partial(\mathcal{CS})$ of the classical $T$-Schottky space generically do not arise from further pinching. Indeed, we show that there are only four values of $\lambda$ for which the $G_\lambda \in \partial(\mathcal{CS})$ correspond to additional cusps.

**Theorem 6.21.** *Let $G_\lambda = \langle S, T \rangle$ be the marked classical $T$-Schottky group generated by parabolic elements $S$ and $T = T_\lambda$. The only words that can be pinched to give an additional parabolic are, up to conjugacy, $S^{\pm 1}T$ and $[S,T]$. At most one of these can be pinched.*

**Proof.** Recall that after replacing $T$ by $T^{-1}$, if necessary, we may assume without loss of generality that $\lambda = |\lambda|e^{i\omega}$, with $0 \leq \omega < \pi$. Let $F$ be the classical $T$-Schottky domain for $G = G_\lambda$, given in Theorem 6.17 by the triple $(\psi, \phi, t)$, where $t = |\tau|$ and

$$\psi = \pi/2 \qquad \phi = \omega/2 - \pi/4 \qquad \tau = e^{i\phi}/\cos\phi.$$

The points of tangency of $F$, not fixed by $S$ or $T$ are as seen in Figure 17:

$$\text{①} \quad \tau(-1-i) \qquad \text{②} \quad \tau(1+i)$$
$$\text{③} \quad \tau(-1+i) \qquad \text{④} \quad \tau(1-i).$$

If a word $W(S,T)$ not conjugate to a power of $S$ or $T$ is parabolic, then by Lemma 6.2, its fixed point must be represented on the boundary of any classical $T$-Schottky domain for $G$ and in particular, on $F$. For this to happen, some of the points of tangency must be paired by $T$.

We consider the various possibilities (see Figure 17).

- $T(\text{①}) = \text{②}$. In this case, $\lambda = \tau(1+i)$, so $\omega = 0$ and hence $\lambda = 2$ and $\tau = 1 - i$. Thus, $S^{-1}T$ is parabolic fixing ①.

- $T(\text{①}) = \text{③}$. In this case, $\lambda = \tau i$, so $\omega = \pi/2$ and hence $\lambda = i$ and $\tau = 1$. It follows that the four points are $G$-equivalent and that $[S,T] = STS^{-1}T^{-1}$ is parabolic fixing ②.

- $T(\text{④}) = \text{②}$. In this case, we again have $\omega = \pi/2$ and hence $\lambda = i$ and $\tau = 1$.

- $T(\text{\textcircled{4}}) = \text{\textcircled{3}}$. In this case, $\omega = \pi$, $\lambda = -2$ and $\tau = 1 + i$.

$\square$

Since the effect of replacing $T$ by $T^{-1}$ is to replace $\lambda$ by $-\lambda$, we have shown

**Corollary 6.22.** *The only values of $\lambda$ on the boundary of classical $T$-Schottky space that arise from additional pinching are $\lambda = \pm i, \pm 2$.*

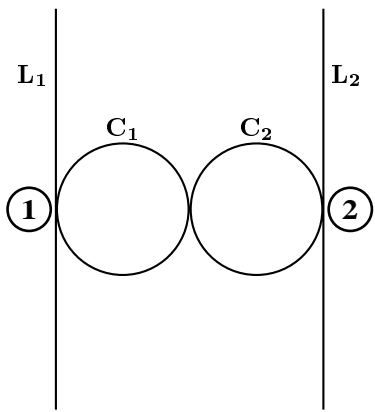
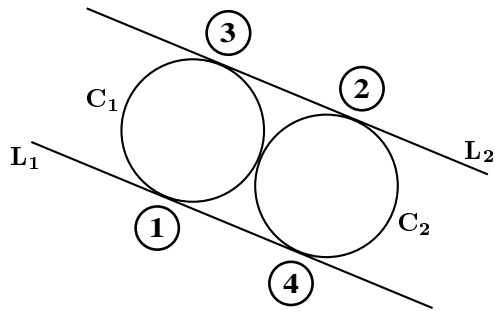

Figure 16: Four tangency extreme classical $T$-Schottky configuration: Classical Domain

Figure 17: Six tangency extreme $T$-Schottky configuration: Classical Domain

**Remark 6.23.** Note that when $\lambda = 2$, the points $\text{\textcircled{3}}$ and $\text{\textcircled{4}}$ are not equivalent. Indeed, $G$ is the ideal triangle group and has a classical $T$-Schottky domain with only four points of tangency, as shown in Figure 16. In this case, the curve $P_3$, shown in Figure 4, has also been pinched; and the noded surface comprises two triply punctured spheres.

**Remark 6.24.** When $\lambda = i$, the points $\text{\textcircled{1}}$, $\text{\textcircled{2}}$, $\text{\textcircled{3}}$ and $\text{\textcircled{4}}$ are all equivalent. In this case, the curve $P_4$, shown in Figure 4, has also been pinched; and the noded surface again comprises two triply punctured spheres.

**Remark 6.25.** Since any classical $T$-Schottky $G_\lambda$ possesses a six tangency domain, a two-parabolic generator Möbius group is classical $T$-Schottky exactly when, up to conjugacy, it has $T$-Schottky domain

$$\{z : |\operatorname{Re}(z)| < 1 \ \& \ |z - i| > 1 \ \& \ |z + i| > 1\}.$$

# 7 Non-Classical $T$-Schottky Groups

We continue with $S$, $T$, $\lambda$ and $\phi$, $\omega$ as before and obtain criteria for the infinitely generated group $\Gamma = \langle T^n S T^{-n} : n \in \mathbb{Z} \rangle$ to be classical $T$-Schottky with equal sized circles. Since $G = \langle S, T \rangle$ is a normal extension of $\Gamma$, the two groups have the same limit set and thus must be simultaneously discrete or non-discrete. However, we show that $\Gamma$ is classical $T$-Schottky at times when $G$ is non-classical $T$-Schottky.

Following [17], we define the region $K$ to be the convex hull of the set consisting of the circle $|z| = 1$ and the points $z = \pm 2$. We prove that $\lambda$ not being in the interior of $K$ is a sufficient condition for $G_\lambda$ to be $T$-Schottky. An immediate consequence is that such groups are free, as was shown in [17].

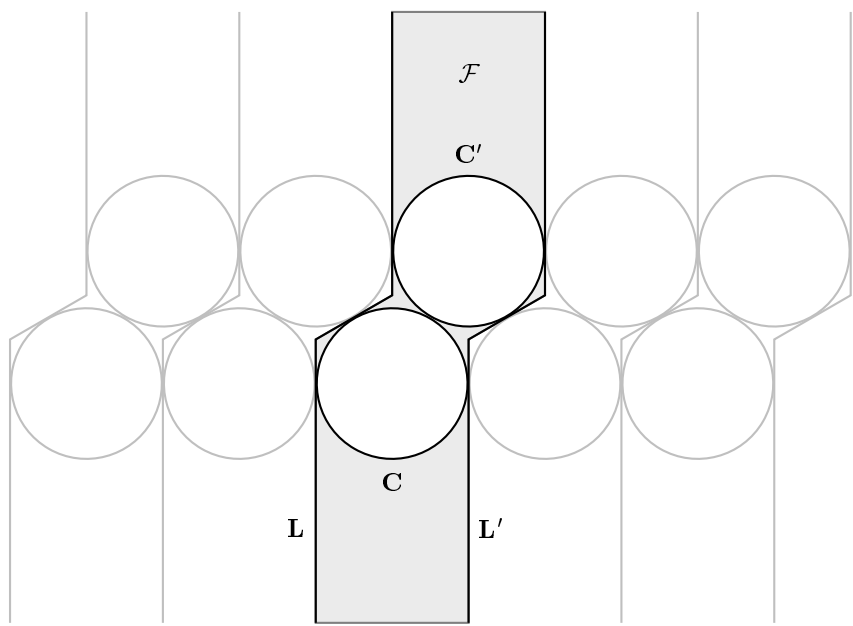

Figure 18: **An Example of a Non-Classical $T$-Schottky Domain.** For the non-classical $T$-Schottky domain $\mathcal{F}$ shown, $\lambda = e^{(i\pi)/3}$, so we normalize via $PSL(2, \mathbb{C})$ to obtain $S = \begin{pmatrix} 1 & 0 \\ e^{(-i\pi/3)} & 1 \end{pmatrix}$ and $T = \begin{pmatrix} 1 & 2 \\ 0 & 1 \end{pmatrix}$. Thus $S$ pairs the circle $C : |z + e^{\frac{i\pi}{3}}| = 1$, with the circle $C' : |z - e^{i\pi/3}| = 1$, and $T$ pairs $L$ with $L'$.

**Lemma 7.1.** *Assume $\Gamma = \langle T^n S T^{-n} : n \in \mathbb{Z} \rangle$ with $0 \leq \omega \leq \pi/2$ is marked classical $T$-Schottky with equal sized circles. Then $\Gamma$ is marked classical $T$-Schottky if and only if*

$$\begin{cases} |\lambda| \geq 1 & \text{and} \quad \pi/3 \leq \omega \leq \pi/2 \\ & \text{or} \\ \operatorname{Re}(\lambda e^{-\pi/3}) \geq 1 & \text{and} \quad 0 \leq \omega \leq \pi/3. \end{cases}$$

**Proof.** $\Gamma$ is a marked classical $T$-Schottky group with equal sized circles if and only if, for some $\phi$, the circles

$$\{|z \pm (e^{i\phi}/\cos\phi) + 2\lambda n| = 1/\cos\phi\}$$

for $n \in \mathbb{Z}$ are disjoint. This holds if and only if

$$|\lambda| \geq 1/\cos\phi \quad \text{and} \quad |\lambda \pm (e^{i\phi}/\cos\phi)| \geq 1/\cos\phi.$$

For a given $\omega$, we therefore need to find the minimum $|\phi|$ for which the following hold:

$$|\lambda| = 1/\cos\phi \quad \text{and} \quad |e^{i\omega} \pm e^{i\phi}| \geq 1.$$

Thus

$$\phi = \begin{cases} 0 & \text{when} \quad \pi/3 \leq \omega \leq \pi/2 \\ \omega - (\pi/3) & \text{when} \quad 0 \leq \omega \leq /\pi/3. \end{cases}$$

The lemma follows. $\square$

Lemma 7.1 gives conditions as to when the infinitely generated group is "marked classical $T$-Schottky." We want to produce a set of Jordan curves for the two generator group. Let $C$ and $D$ be the circles of equal radii passing through the origin, with $S$ mapping the exterior of $C$ to the interior of $D$. The inequalities assure that there is a Jordan curve $\gamma$ passing through $\infty$, separating the interiors of $C$ and $D$ from those of $T(C)$ and $T(D)$, and such that the curves $C$, $D$, $T^{-1}(\gamma)$, and $\gamma$ bound a non-classical $T$-Schottky domain for $G$. Note that $\gamma$ may be tangent to one or all of $C$, $D$, $T(C)$ and $T(D)$. This shows that the two generator group is non-classical $T$-Schottky.

An example of an explicit non-classical $T$-Schottky domain is shown in Figure 18.

In Figure 15, the set of classical $T$-Schottky groups is shown with its bounding parabolas superimposed on the Riley Slice. One can thus view (numerically) the set of non-classical $T$-Schottky groups.

As a consequence of Lemma 7.1, we have a precise description of a subset of these shown in dark gray in Figure 7.

Let K be the set considered by Lyndon and Ullman in Figure 5. Namely, if $\lambda$ is in the first quadrant it must satisfy

$$\begin{cases} |\lambda| < 1 & \text{and} \quad \pi/3 \leq \omega \leq \pi/2 \\ & \text{or} \\ \operatorname{Re}(\lambda e^{-\pi/3}) < 1 & \text{and} \quad 0 \leq \omega \leq \pi/3. \end{cases}$$

**Theorem 7.2.** *If $\lambda$ lies in the first quadrant below the $T$-Schottky parabola and exterior to $K$, then $G_\lambda$ is a non-classical $T$-Schottky group.*

The above may be phrased as a criterion for a Möbius group $G$ to be $T$-Schottky.

**Theorem 7.3.** *If $G = \langle S, T \rangle$ with $tr(S) = tr(T) = 2$, then $G$ is $T$-Schottky if*

$$\begin{cases} |tr(ST) - 2| \geq 2 & \text{and} \quad \pi/3 \leq \operatorname{Arg}[tr(ST) - 2)] \leq \pi/2, \\ & \text{or} \\ \operatorname{Re}[(tr(ST) - 2)e^{-\pi/3}] \geq 2 & \text{and} \quad 0 \leq \operatorname{Arg}[tr(ST) - 2] \leq \pi/3. \end{cases}$$

**Proof.** The result follows by recalling that $\lambda = (1/2)[tr(ST) - 2]$ □

# 8 Marked groups with Loxodromic - Parabolic generators.

In order to complete the description of the marked classical $T$-Schottky spaces, we need to describe the boundary for case (ii) of Theorem 6.5.

We define $R = \begin{pmatrix} 1 & 2\lambda \\ -1 & (1-2\lambda) \end{pmatrix}$ and $\Gamma_\lambda = \langle R, T \rangle$ and note that $TR^{-1} = \begin{pmatrix} 1 & 0 \\ 1 & 1 \end{pmatrix} = S$.

Hence, as un-marked groups $G_\lambda = \Gamma_\lambda$. Further, the following are equivalent:

- $\Gamma_\lambda$ is marked classical $T$-Schottky
- $\Gamma_{\bar{\lambda}}$ is marked classical $T$-Schottky
- $\Gamma_{-\lambda}$ is marked classical $T$-Schottky
- $\Gamma_{-\bar{\lambda}}$ is marked classical $T$-Schottky.

As depicted in Figure 19, we can now describe when $\Gamma_\lambda$ is marked classical.

**Theorem 8.1.** *$\Gamma_\lambda$ is marked classical $T$-Schottky if and only if $|Re(\lambda)| + |Im(\lambda)| \geq 2$.*

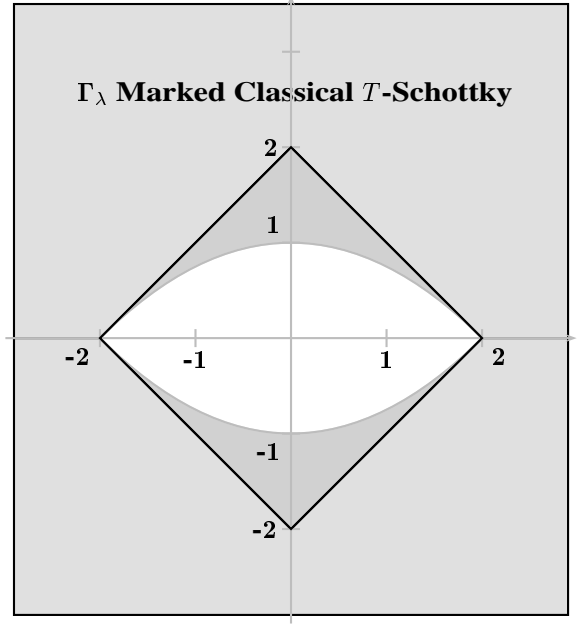

Figure 19: **Marked Classical $T$-Schottky Space.** Depicted is the complex $\lambda$-plane. For $\lambda$ in the dark gray area $\Gamma_\lambda$ is Classical $T$-Schottky but not marked Classical $T$-Schottky on the generators $R$ and $T$.

**Proof.** By the above, we may assume that $\lambda$ lies in the first quadrant.

Any $T$-Schottky domain $\mathcal{F}$ for $\Gamma_\lambda$ is bounded by lines $L_1$ and $L_2 = T(L_1)$ together with circles $C_1$ and $C_2 = R(C_1)$, as shown in Figure 20. Since $S(0) = 0$, $L_2$ must be tangent to $C_2$ at the parabolic fixed point $0$ by Lemma 6.2. Similarly, $C_1$ must be tangent to $L_1$ at $T^{-1}(0) = -2\lambda$.

We thereby obtain from Corollary 6.9 that for some $\tau$,

$$C_1 = \left\{ z : \left| z - \left( \frac{\tau}{\tau + \bar{\tau} - 1} - 2\lambda \right) \right| = \frac{|\tau|}{\tau + \bar{\tau} - 1} text with \tau + \bar{\tau} - 1 > 0 \right\}$$

and

$$C_2 = \{ z : |z + \tau| = |\tau| \}.$$

$\mathcal{F}$ is then a $T$-Schottky domain for the marked group $\Gamma_\lambda$ precisely when both of the following conditions hold:

- $C_1$ does not overlap $C_2$, i.e., if

$$\left| \lambda - \frac{(\tau + \bar{\tau})\tau}{2(\tau + \bar{\tau} - 1)} \right| \geq \left| \frac{(\tau + \bar{\tau})\tau}{2(\tau + \bar{\tau} - 1)} \right|;$$

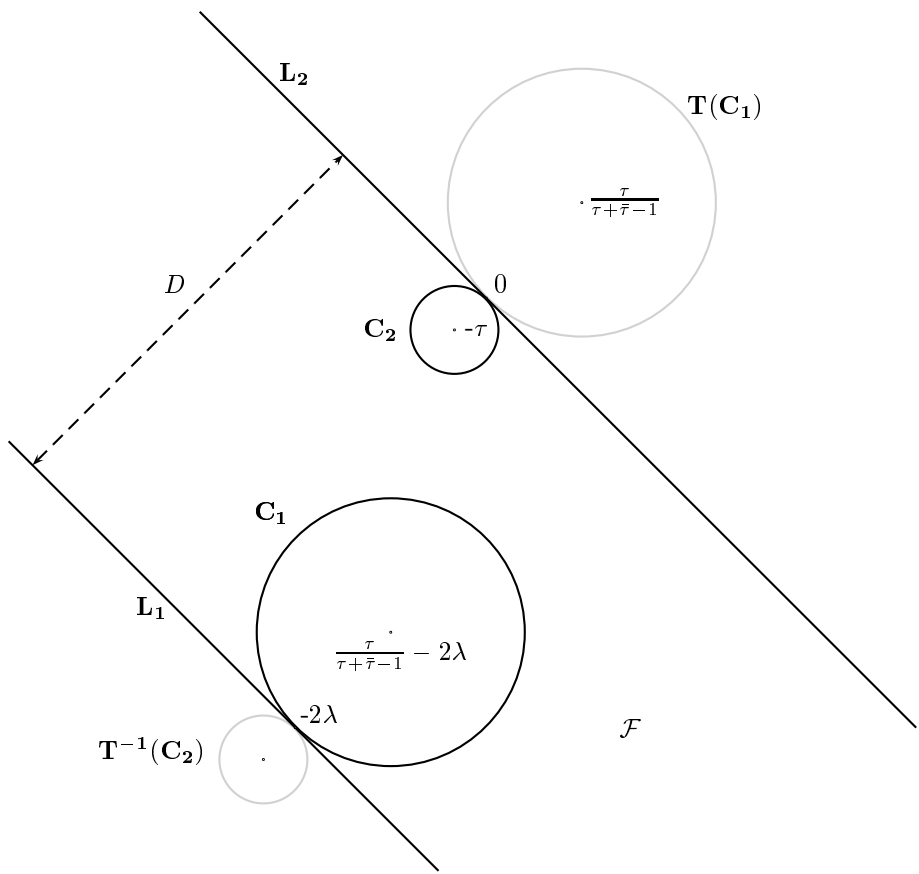

Figure 20: $T$-Schottky Domain: Loxodromic - Parabolic generators

- $C_1$ and $C_2$ lie between $L_1$ and $L_2$, i.e., if

$$(\lambda - m\tau) \cdot \tau \geq 0,$$

where $m = \max\{1, 1/(\tau + \bar{\tau} - 1)\}$ and $\cdot$ denotes the dot product.

For fixed $\tau$, the possible $\lambda$ are shown in Figure 21.

If the argument of $\tau$ is fixed, then by convexity choosing $\tau + \bar{\tau} - 1 = 1$ achieves the largest domain of possible $\lambda$. Observe that this choice is equivalent to requiring that $C_1$ and $C_2$ have the same radius.

For fixed $\lambda$, $\Gamma_\lambda$ thus possesses a $T$-Schottky domain $\mathcal{F}$ precisely when there exists $\tau$ with $Re(\tau) = 1$ such that

- $C_1$ does not overlap $C_2$, i.e.,

$$|\tau - \lambda| \geq |\tau|;$$

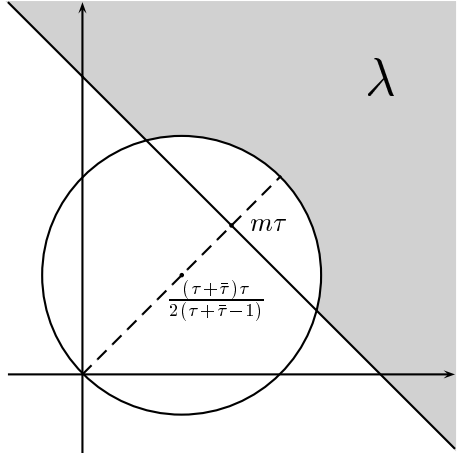
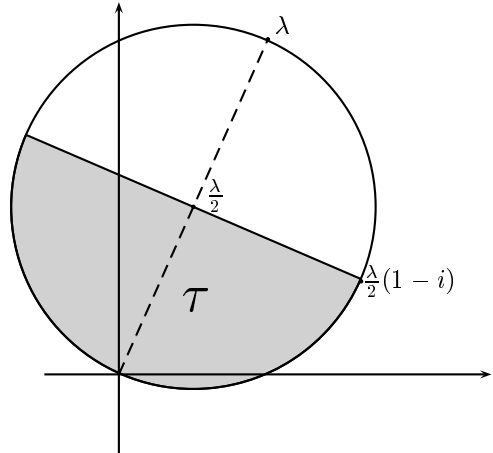

Figure 21: $\lambda$ - Separation        Figure 22: $\tau$ - Separation

- $C_1$ and $C_2$ lie between $L_1$ and $L_2$, i.e.,

$$|\tau - \lambda/2| \leq |\lambda/2|.$$

The above pair of conditions are shown in Figure 22.

It follows that $\Gamma_\lambda$ thus possesses a $T$-Schottky domain $\mathcal{F}$ precisely when $Re\left(\frac{\lambda}{2}(1-i)\right) \geq 1$. If this holds, then $\tau = 1 + iIm\left(\frac{\lambda}{2}(1-i)\right)$ suffices. The theorem follows. □

**8.1 Technical lemmas**   The following are elementary computations.

**Lemma 8.2.** *Let $a_0$, $b_0$ and $w$ be complex numbers and $R$ a real number with $R \neq 1$. Then*

$$|w + a_0| = R|w + b_0|$$
$$\Leftrightarrow \left|w - \frac{a_0 - R^2 b_0}{R^2 - 1}\right| = \frac{R|a_0 - b_0|}{|R^2 - 1|}.$$

**Lemma 8.3.**

$S = \begin{pmatrix} \nu & 0 \\ \nu^{-1} & \nu^{-1} \end{pmatrix}$ *maps the circle* $C_1 : |z - \alpha| = r$ *to the circle* $C_2 : |z - \beta| = \rho$

*if and only if*

$$\beta = \frac{\nu^2(r^2 - \alpha - |\alpha|^2)}{r^2 - |1 + \alpha|^2} \quad and \quad \rho = \frac{r|\nu|^2}{|r^2 - |1 + \alpha|^2|}.$$

*Further, $S$ maps the exterior of $C_1$ to the interior of $C_2$ $\iff$ $|\alpha + 1| < r$.*

# Acknowledgment

The authors wish to thank T.Jørgensen for many helpful conversations during the preparation of this work. We thank D. Wright for giving us permission to use his drawing of the Riley slice.

*Jane Gilman*
???????????
  ?????????,
    ?????????
      email: gilman@andromeda.rutgers.edu

*Peter Waterman*
?????????
  ?????????
    ?????????
      email: waterman@math.niu.edu